\documentclass[letterpaper, 11pt]{article}
\usepackage{amsmath,amsfonts,amssymb,amsthm,xypic,xy,hyperref}
\xyoption{all}

\newtheorem{theorem}{Theorem}[section]
\newtheorem{corollary}[theorem]{Corollary}

\newtheorem{prop}[theorem]{Proposition}

\newtheorem*{remark}{Remark}
\theoremstyle{definition}
\newtheorem{defn}[theorem]{Definition}
\newtheorem{example}[theorem]{Example}
\numberwithin{equation}{section}

\newcommand{\nc}{\newcommand}\nc{\br}{\overline}
\nc{\FH}{\mathcal H} \nc{\CC}{\mathbb C}\nc{\DD}{\mathcal D}
\nc{\Cc}{\mathcal C}\nc{\JJ}{\mathcal J}\nc{\II}{\mathcal I}
\nc{\PP}{\mathbb P} \nc{\KK}{\mathbb K} \nc{\RR}{\mathbb R}
\nc{\LL}{\mathcal L} \nc{\Ll}{\ell} \nc{\NN}{\mathbb
N}\nc{\TT}{\mathcal{T}} \nc{\ZZ}{\mathbb Z} \nc {\HH}{\mathbb H} \nc
{\OO}{\mathcal{O}} \nc{\lra}{\longrightarrow} \nc{\bdot}{\bullet}
\nc{\w}{\omega}
\nc{\dd}{\mathcal{D}}\nc{\im}{\mathrm{Im}}\nc{\Nij}{\mathrm{Nij}}\nc{\Jac}{\mathrm{Jac}}

\newcommand{\rk}{\mathrm{rk\ }}\nc{\sgn}{\mathrm{sgn}}
\newcommand{\Tr}{\mathrm{Tr}}
\newcommand{\End}{\mathrm{End}}
\newcommand{\Ann}{\mathrm{Ann}}
\newcommand{\delbar}{\overline{\partial}}
\newcommand{\so}{\mathfrak{so}}\nc{\MM}{\mathcal{M}}\nc{\Har}{\mathcal{H}} \nc{\zed}{\mathcal{Z}}
\nc{\image}{\mathrm{Im}\ }
\newcommand{\IP}[1]{\langle #1\rangle}

\nc{\ba}{\overline} \nc{\del}{\partial} \nc{\AAA}{\mathcal{A}}
\nc{\de}{\delta}\nc{\debar}{\overline{\delta}} \nc{\FF}{\mathcal{F}}
\nc{\Ja}{e^{\tfrac{\pi}{2}\JJ_1}} \nc{\Jb}{e^{\tfrac{\pi}{2}\JJ_2}}
\nc{\eps}{\epsilon}\nc{\id}{\mathrm{id}}\nc{\Dir}{\mathrm{Dir}}\nc{\SO}{\mathrm{SO}}
\nc{\IPS}[1]{#1}\nc{\veps}{\varepsilon}\nc{\Cour}[1]{[#1]}\nc{\Diff}{\mathrm{Diff}}
\nc{\ad}{\mathrm{ad}} \nc{\UU}{\mathcal{U}}\nc{\Aa}{\mathcal
A}\nc{\Pic}{\mathrm{Pic}}\nc{\type}{\mathrm{type}}
\title{\bf Generalized complex geometry}
\author{Marco Gualtieri}
\date{}
\begin{document}
\maketitle \abstract{Generalized complex geometry encompasses
complex and symplectic geometry as its extremal special cases.  We
explore the basic properties of this geometry, including its
enhanced symmetry group, elliptic deformation theory, relation to
Poisson geometry, and local structure theory.  We also define and
study generalized complex branes, which interpolate between flat
bundles on Lagrangian submanifolds and holomorphic bundles on
complex submanifolds.}

\tableofcontents
\section*{Introduction}

Generalized complex geometry arose from the work of
Hitchin~\cite{Hitchin} on geometries defined by stable differential
forms of mixed degree.  Algebraically, it interpolates between a
symplectic form $\omega$ and a complex structure $J$ by viewing each
as a complex (or equivalently, symplectic) structure $\JJ$ on the
direct sum of the tangent and cotangent bundles $T\oplus T^*$,
compatible with the natural split-signature metric which exists on
this bundle.  Remarkably, there is an integrability condition on
such generalized complex structures which specializes to the closure
of the symplectic form on one hand, and the vanishing of the
Nijenhuis tensor of $J$ on the other. This is simply that $\JJ$ must
be integrable with respect to the \emph{Courant bracket}, an
extension of the Lie bracket of vector fields to smooth sections of
$T\oplus T^*$ which was introduced by Courant and
Weinstein~\cite{Courant},\cite{CourWein} in their study of Dirac
structures.  Dirac structures interpolate between Poisson bivector
fields and closed 2-forms; in this sense, generalized complex
geometry is a complex analogue of Dirac geometry.

We begin, in parts~\ref{linalg} and~\ref{courant}, with a study of
the natural split-signature orthogonal structure on $T\oplus T^*$
and its associated spin bundle $\wedge^\bullet T^*$, the bundle of
differential forms.  Viewing forms as spinors then leads to a
natural definition of the Courant bracket, and we study two
remarkable properties of this bracket; first, its symmetry group is
an extension of the group of diffeomorphisms by the abelian group of
closed 2-forms (B-field transformations), and second, it can be
``twisted'' by a real closed 3-form $H$.  We describe what this
means in the language of gerbes and exact Courant algebroids.  We
also provide a brief review of Dirac geometry and introduce the
notion of tensor product of Dirac structures, obtained independently
by Alekseev-Bursztyn-Meinrenken in~\cite{ABM}.

In part~\ref{gcg}, we treat the basic properties of generalized
complex structures. We show that any generalized complex manifold
admits almost complex structures, and has two natural sets of Chern
classes $c_k^\pm$. We also show that a generalized complex structure
is determined by a complex pure spinor line subbundle $K\subset
\wedge^\bullet T^*\otimes\CC$, the \emph{canonical bundle}, which
can be seen as the minimal degree component of an induced
$\ZZ$-grading on the (usually $\ZZ/2\ZZ$-graded) twisted de Rham
complex $(\Omega^\bullet(M), d_H)$, where $d_H = d + H\wedge\cdot$.
For a complex structure, $K$ is the usual canonical bundle, wheras
for a symplectic structure it is generated by the mixed-degree form
$e^{i\omega}$.  We also describe a real Poisson structure $P$
naturally associated to any generalized complex structure and
discuss its modular class.  We conclude with an example of a family
of generalized complex structures which interpolates between the
complex and symplectic structure on a holomorphic symplectic
manifold.

In part~\ref{darboux}, we prove a local structure theorem for
generalized complex manifolds, analogous to the Darboux theorem in
symplectic geometry and the Newlander-Nirenberg theorem in complex
geometry. We show that near any regular point for the Poisson
structure $P$, the generalized complex manifold is equivalent, via a
diffeomorphism and a B-field transformation, to a product of a
complex space of dimension $k$ (called the \emph{type}) with a
symplectic space. Finally, we provide an example of a generalized
complex manifold whose type is constant outside of a submanifold,
along which it jumps to a higher value.

In part~\ref{deform}, we develop the deformation theory of
generalized complex manifolds.  It is governed by a differential
Gerstenhaber algebra $(C^\infty(\wedge^k L^*),d_L, [\cdot,\cdot])$
constructed from the $+i$-eigenbundle $L$ of $\JJ$. This
differential complex is elliptic, and therefore has
finite-dimensional cohomology groups $H^k(M,L)$ over a compact
manifold $M$.  Integrable deformations correspond to sections
$\eps\in C^\infty(\wedge^2 L^*)$ satisfying the Maurer-Cartan
equation
\begin{equation*}
d_L\eps+\tfrac{1}{2}[\eps,\eps]=0.
\end{equation*}
Similarly to the case of deformations of complex structure, there is
an analytic obstruction map $\Phi:H^2(M,L)\rightarrow H^3(M,L)$, and
if this vanishes then there is a locally complete family of
deformations parametrized by an open set in $H^2(M,L)$. In the case
that we are deforming a complex structure, this cohomology group
turns out to be
\begin{equation*}
H^0(M,\wedge^2\mathcal{T})\oplus H^1(M,\mathcal{T})\oplus
H^2(M,\mathcal{O}).
\end{equation*}
This is familiar as the ``extended deformation space'' of Barannikov
and Kontsevich~\cite{KontsevichBarannikov}, for which a geometrical
interpretation has been sought for some time.  Here it appears
naturally as the deformation space of a complex structure in the
generalized setting.

Finally, in part~\ref{branes}, we introduce \emph{generalized
complex branes}, which are vector bundles supported on submanifolds
for which the pullback of the ambient gerbe is trivializable,
together with a natural requirement of compatibility with the
generalized complex structure.  The definition is similar to that of
a D-brane in physics; indeed, we show that for a usual symplectic
manifold, generalized complex branes consist not only of flat vector
bundles supported on Lagrangian submanifolds, but also certain
bundles over a class of coisotropic submanifolds equipped with a
holomorphic symplectic structure transverse to their characteristic
foliation. These are precisely the co-isotropic A-branes discovered
by Kapustin and Orlov~\cite{Kapustin}.

This article is largely based upon the author's doctoral
thesis~\cite{Gu}, supported by the Rhodes Trust.  Thanks are due
especially to Nigel Hitchin for his guidance and insight. Many
thanks as well to Henrique Bursztyn, Gil Cavalcanti, Jacques
Hurtubise, Anton Kapustin, and Alan Weinstein for helpful
conversations.

\section{Linear geometry of $T\oplus T^*$}\label{linalg}

Let $V$ be a real vector space of dimension $m$, and let $V^*$ be
its dual space.  Then $V\oplus V^*$ is endowed with a natural
symmetric bilinear form
\begin{align*}
\IP{X+\xi,Y+\eta}&=\tfrac{1}{2}(\xi(Y)+\eta(X)).
\end{align*}
This non-degenerate inner product has signature $(m,m)$ and
therefore has symmetry group $\mathrm{O}(V\oplus V^*)\cong
\mathrm{O}(m,m)$, a non-compact orthogonal group. In addition,
$V\oplus V^*$ has a canonical orientation, since $\det(V\oplus V^*)
= \det V\otimes\det V^*=\RR$.

\subsection{Symmetries of $V\oplus V^*$}
The Lie algebra of the special orthogonal group $SO(V\oplus V^*)$ is
defined as usual:
\begin{align*}
\mathfrak{so}(V\oplus V^*)=\left\{R\ \big|\ \IP{Tx,y}+\IP{x,Ty}=0\ \
\forall\  x,y\in V\oplus V^*\right\}.
\end{align*}
Using the splitting $V\oplus V^*$ we can decompose as follows:
\begin{equation}\label{comp}
T=\left(\begin{matrix}A&\beta\\B&-A^*\end{matrix}\right),
\end{equation}
where $A\in \End(V)$, $B:V\lra V^*$, $\beta:V^*\lra V$, and where
$B$ and $\beta$ are skew, i.e. $B^*=-B$ and $\beta^*=-\beta$.
Therefore we may view $B$ as a 2-form in $\wedge^2V^*$ via $B(X)=i_X
B$ and similarly we may regard $\beta$ as an element of $\wedge^2
V$, i.e. a bivector.  This corresponds to the observation that
$\so(V\oplus V^*)=\wedge^2(V\oplus V^*)=\End(V)\oplus \wedge^2
V^*\oplus\wedge^2V$.

By exponentiation, we obtain certain important orthogonal symmetries
of $T\oplus T^*$ in the identity component of $SO(V\oplus V^*)$.
\begin{example}[$B$-transform]
First let $B$ be as above, and consider
\begin{equation}
\exp(B)=\left(\begin{matrix}1&\\B&1\end{matrix}\right),
\end{equation}
an orthogonal transformation sending $X+\xi\mapsto X+\xi+i_XB$. It
is useful to think of $\exp(B)$ a \emph{shear} transformation, which
fixes projections to $V$ and acts by shearing in the $V^*$
direction.
\end{example}
\begin{example}[$\beta$-transform]\label{betatrans}
Similarly, let $\beta$ be as above, and consider
\begin{equation}
\exp(\beta)=\left(\begin{matrix}1&\beta\\&1\end{matrix}\right),
\end{equation}
an orthogonal transformation sending $X+\xi\mapsto
X+\xi+i_\xi\beta$. Again, it is useful to think of $\exp(\beta)$ a
shear transformation, which fixes projections to $T^*$ and acts by
shearing in the $T$ direction.
\end{example}
\begin{example}
If we choose $A\in \End(V)$ as above, then since
\begin{equation}
\exp(A)=\left(\begin{matrix}\exp A&\\&(\exp
A^*)^{-1}\end{matrix}\right),
\end{equation}
we have a distinguished embedding of $GL^+(V)$ into the identity
component of $SO(V\oplus V^*)$.
\end{example}

\subsection{Maximal isotropic subspaces}\label{maxiso}

A subspace $L\subset V\oplus V^*$ is \emph{isotropic} when
$\IP{x,y}=0$ for all $x,y\in L$. Since we are in signature $(m,m)$,
the maximal dimension of such a subspace is $m$, and if this is the
case, $L$ is called maximal isotropic. Maximal isotropic subspaces
of $V\oplus V^*$ are also called \emph{linear Dirac structures}
(see~\cite{Courant}).  Note that $V$ and $V^*$ are examples of
maximal isotropics. The space of maximal isotropics has two
connected components, and elements of these are said to have
\emph{even} or \emph{odd} parity, depending on whether they share
their connected component with $V$ or not, respectively. This
situation becomes more transparent after studying the following two
examples.
\begin{example}
Let $\Delta\subset V$ be any subspace.  Then the subspace
\begin{equation*}
\Delta\oplus\Ann(\Delta) \subset V\oplus V^*,
\end{equation*}
where $\Ann (\Delta)$ is the annihilator of $\Delta$ in $V^*$, is
maximally isotropic.
\end{example}
\begin{example}
Let $i:\Delta\hookrightarrow V$ be a subspace inclusion, and let
$\veps\in \wedge^2 \Delta^*$. Then the subspace
\begin{equation*}
L(\Delta,\veps)=\left\{X+\xi\in \Delta\oplus V^*\ :\
i^*\xi=i_X\veps\right\}\subset V\oplus V^*
\end{equation*}
is an extension of the form
\begin{equation*}
\xymatrix{0\ar[r]&\Ann(\Delta)\ar[r]&L(\Delta,\veps)\ar[r]&\Delta\ar[r]&0},
\end{equation*}
and satisfies
$\IP{X+\xi,Y+\eta}=\tfrac{1}{2}(\veps(Y,X)+\veps(X,Y))=0$ for all
$X+\xi,Y+\eta\in L(\Delta,\veps)$, showing that $L(\Delta,\veps)$ is
a maximal isotropic subspace.
\end{example}
Note that the second example specializes to the first by taking
$\eps=0$. Furthermore note that $L(V,0)=V$ and $L(\{0\},0)=V^*$. It
is not difficult to see that every maximal isotropic is of this
form:
\begin{prop}\label{Eeps}
Every maximal isotropic in $V\oplus V^*$ is of the form
$L(\Delta,\veps)$.
\end{prop}
\begin{proof}
Let $L$ be a maximal isotropic and define $\Delta=\pi_VL$, where
$\pi_V$ is the canonical projection $V\oplus V^*\lra V$. Then since
$L$ is maximal isotropic, $L\cap V^*=\Ann(\Delta)$.  Finally note
that $\Delta^*=V^*/\Ann(\Delta)$, and define
$\veps:\Delta\rightarrow \Delta^*$ via $\veps: x\mapsto
\pi_{V^*}(\pi_V^{-1}(x)\cap L)\in V^*/\Ann(\Delta)$. Then
$L=L(\Delta,\veps)$.
\end{proof}

The integer $k=\dim\Ann(\Delta)=m-\dim\pi_V(L)$ is a useful
invariant associated to a maximal isotropic in $V\oplus V^*$, and
determines the parity as we now explain.
\begin{defn}\label{type}
The \emph{type} of a maximal isotropic $L\subset V\oplus V^*$ is the
codimension $k$ of its projection onto $V$.
\end{defn}
Since a $B$-transform preserves projections to $V$, it does not
affect $\Delta$:
\begin{align*}
\exp B(L(\Delta,\veps))&=L(\Delta,\veps+i^*B),
\end{align*}
where $i:\Delta\hookrightarrow V$ is the inclusion.  Hence
$B$-transforms do not change the type of the maximal isotropic. In
fact, we see that by choosing $B$ and $\Delta$ accordingly, we can
obtain any maximal isotropic as a $B$-transform of $L(\Delta,0)$.

On the other hand, $\beta$-transforms do modify projections to $V$,
and therefore may change the dimension of $\Delta$.  To see how this
occurs, we express the maximal isotropic as a generalized graph from
$V^*\rightarrow V$, i.e. define $F=\pi_{V^*}L$ and
$\gamma\in\wedge^2F^*$ given by
$\gamma(y)=\pi_{V}(\pi_{V^*}^{-1}(y)\cap L)$ modulo $\Ann(F)$. As
before, define
\begin{equation}\label{othproj}
L(F,\gamma)=\left\{X+\xi\in V\oplus F\ :\ j^*X=i_\xi\gamma\right\},
\end{equation}
where $j:F\hookrightarrow V^*$ is the inclusion.  Now, the
projection $\Delta=\pi_VL(F,\gamma)$ always contains $L\cap
V=Ann(F)$, and if we take the quotient of $\Delta$ by this subspace
we obtain the image of $\gamma$ in $F^*=V/\Ann(F)$:
\begin{equation*}
\frac{\Delta}{L\cap V}=\frac{\Delta}{\Ann(F)}=\mathrm{Im}(\gamma).
\end{equation*}
Therefore, we obtain the dimension of $\Delta$ as a function of
$\gamma$:
\begin{equation*}
\dim \Delta=\dim L\cap V+\rk\gamma.
\end{equation*}
Because $\gamma$ is a skew form, its rank is even. A
$\beta$-transform sends $\gamma\mapsto\gamma+j^*\beta$, which also
has even rank, and therefore we see that a $\beta$-transform, which
is in the identity component of $SO(V\oplus V^*)$, can be used to
change the dimension of $\Delta$, and hence the type of
$L(\Delta,\eps)$, by an even number, yielding the following result.
\begin{prop}\label{partype} The space of maximal isotropics in $V\oplus V^*$ has
two connected components, distinguished by the parity of the type;
maximal isotropics have even parity if and only if they share a
connected component with $V$.
\end{prop}

\subsection{Spinors for $V\oplus V^*$: exterior
forms}\label{spinandforms}

Let $CL(V\oplus V^*)$ be the Clifford algebra defined by the
relation
\begin{equation}\label{clifford}
v^2=\IP{v,v},\ \ \forall v\in V\oplus V^*,
\end{equation}
The Clifford algebra has a natural representation on
$S=\wedge^\bullet V^*$ given by
\begin{equation}\label{cliff}
(X+\xi)\cdot\varphi = i_X\varphi + \xi\wedge\varphi,
\end{equation}
where $X+\xi\in V\oplus V^*$ and $\varphi\in\wedge^\bullet V^*$. We
verify that this defines an algebra representation:
\begin{align*}
(X+\xi)^2\cdot\varphi &=
i_X(i_X\varphi+\xi\wedge\varphi)+\xi\wedge(i_X\varphi+\xi\wedge\varphi)\\
&=(i_X\xi)\varphi\\
&=\IP{X+\xi,X+\xi}\varphi,
\end{align*}
as required.  This representation is the standard spin
representation, so that $CL(V\oplus V^*)=\End(\wedge^\bullet V^*)$.
Since in signature $(m,m)$ the volume element $\omega$ of a Clifford
algebra satisfies $\omega^2=1$, the spin representation decomposes
into the $\pm1$ eigenspaces of $\omega$ (the positive and negative
helicity spinors):
\begin{equation*}
S=S^+\oplus S^-,
\end{equation*}
corresponding to the decomposition
\begin{equation*}
\wedge^\bullet V^*=\wedge^\text{ev}V^*\oplus\wedge^\text{od}V^*
\end{equation*}
according to parity. While the splitting $S=S^+\oplus S^-$ is not
preserved by the whole Clifford algebra, $S^\pm$ are irreducible
representations of the spin group, which sits in the Clifford
algebra as
\begin{equation*}
\mathrm{Spin}(V\oplus V^*)=\{v_1\cdots v_r\ :\ v_i\in V\oplus V^*,
\IP{v_i,v_i}=\pm 1\ \text{and}\ r\ \text{even}\},
\end{equation*}
and which is a double cover of $SO(V\oplus V^*)$ via the
homomorphism
\begin{align}
\rho&:\mathrm{Spin}(V\oplus V^*)\lra SO(V\oplus V^*)\notag\\
\rho&(x)(v)=xv x^{-1},\ \ x\in \mathrm{Spin}(V\oplus V^*),\ v\in
V\oplus V^*.\label{repn}
\end{align}
We now describe the action of $\so (V\oplus V^*)$ in the spin
representation. Recall that $\mathfrak{so}(V\oplus
V^*)=\wedge^2(V\oplus V^*)$ sits naturally inside the Clifford
algebra, and that the derivative of $\rho$, given by
\begin{equation*}
d\rho_x(v)=xv-vx = [x,v],\ \ x\in \mathfrak{so}(V\oplus V^*),\ \
v\in V\oplus V^*,
\end{equation*}
defines the usual representation of $\mathfrak{so}(V\oplus V^*)$ on
$V\oplus V^*$.  In the following we take $\{e_i\}$ to be a basis for
$V$ and $\{e^i\}$ the dual basis.

\begin{example}[$B$-transform]
A 2-form $B=\frac{1}{2}B_{ij}e^i\wedge e^j$, $B_{ij}=-B_{ji}$, has
image in the Clifford algebra given by $\frac{1}{2}B_{ij}e^je^i$,
and hence its spinorial action on an exterior form
$\varphi\in\wedge^\bullet V^*$ is
\begin{equation*}
B\cdot\varphi =
\frac{1}{2}B_{ij}e^j\wedge(e^i\wedge\varphi)=-B\wedge\varphi.
\end{equation*}
Exponentiating, we obtain
\begin{equation}\label{Bfield}
e^{-B}\varphi=(1-B+\tfrac{1}{2}B\wedge B + \cdots)\wedge\varphi.
\end{equation}
\end{example}
\begin{example}[$\beta$-transform]
A bivector $\beta=\frac{1}{2}\beta^{ij}e_i\wedge e_j$,
$\beta^{ij}=-\beta^{ji}$, has image in the Clifford algebra given by
$\frac{1}{2}\beta^{ij}e^je^i$. Its spinorial action on a form
$\varphi$ is
\begin{equation*}
\beta\cdot\varphi =
\frac{1}{2}\beta^{ij}i_{e_j}(i_{e_i}\varphi)=i_{\beta}\varphi.
\end{equation*}
Exponentiating, we obtain
\begin{equation}
e^\beta\varphi=(1+i_\beta+\tfrac{1}{2}i_\beta^2 + \cdots)\varphi.
\end{equation}
\end{example}

\begin{example}[$GL^+(V)$ action]
An endomorphism $A=A_i^je^i\otimes e_j$ has image in the Clifford
algebra given by $\frac{1}{2}A_i^j(e_je^i-e^ie_j)$, and has
spinorial action
\begin{align*}
A\cdot\varphi&=\tfrac{1}{2}A_i^j(i_{e_j}(e^i\wedge\varphi)-e^i\wedge
i_{
e_j}\varphi)\\
&=\tfrac{1}{2}A_i^j\delta^i_j\varphi - A_i^je^i\wedge
i_{e_j}\varphi\\
&=-A^*\varphi + \tfrac{1}{2}\Tr (A) \varphi,
\end{align*}
where $\varphi\mapsto -A^*\varphi=-A^j_i e^i\wedge i_{e_j}\varphi$
is the usual action of $End(V)$ on $\wedge^\bullet V^*$.  Hence, by
exponentiation, the spinorial action of $GL^+(V)$ on $\wedge^\bullet
V^*$ is by
\begin{equation*}
g\cdot\varphi=\sqrt{\det g}(g^*)^{-1}\varphi,
\end{equation*}
i.e. as a $GL^+(V)$ representation the spinor representation
decomposes as
\begin{equation}\label{spinors}
S=\wedge^\bullet V^*\otimes (\det V)^{1/2}.
\end{equation}
\end{example}
In fact, as we shall see in the following sections, tensoring with
the half densities as in~\eqref{spinors} renders $S$ independent of
polarization, i.e. if $N,N'$ are maximal isotropic subspaces such
that $V\oplus V^* = N+N'$, then the inner product gives an
identification $N'=N^*$, and one obtains a canonical isomorphism
$\wedge^\bullet V^*\otimes (\det V)^{1/2} = \wedge^\bullet
N^*\otimes(\det N)^{1/2}$.

\subsection{The Mukai pairing}

There is an invariant bilinear form on spinors, which we now
describe, following the treatment of Chevalley~\cite{Chevalley}. For
$V\oplus V^*$ this bilinear form coincides with the Mukai pairing of
forms~\cite{Mukai}.

Since we have the splitting $V\oplus V^*$ into maximal isotropics,
the exterior algebras $\wedge^\bullet V$ and $\wedge^\bullet V^*$
are subalgebras of $CL(V\oplus V^*)$.  In particular, $\det V$ is a
distinguished line inside $CL(V\oplus V^*)$, and it generates a left
ideal with the property that upon choosing a generator $f\in\det V$,
every element of the ideal has a unique representation as $sf$,
$s\in\wedge^\bullet V^*$.  This defines an action of the Clifford
algebra on $\wedge^\bullet V^*$ by
\begin{equation*}
(x\cdot s)f=xsf\ \ \ \forall x\in CL(V\oplus V^*),
\end{equation*}
which is the same action as that defined by~(\ref{cliff}).

Let $x\mapsto x^\top$ denote the main antiautomorphism of the
Clifford algebra, i.e. that determined by the tensor map
$v_1\otimes\cdots\otimes v_k\mapsto v_k\otimes \cdots \otimes v_1$.
Now let $s,t\in \wedge^\bullet V^*$ and consider the \emph{Mukai
pairing} $(\cdot,\cdot):\xymatrix{\wedge^\bullet V^*\otimes
\wedge^\bullet V^*\ar[r]&\det V^*}$ given by
\begin{equation*}
\IPS{(s,t)}=[s^\top\wedge t]_{m},
\end{equation*}
where $[\cdot]_{m}$ indicates projection to the component of degree
$m=\dim V$. We can express $(\IPS{,})$ in the following way, using
any generator $f\in\det V$:
\begin{equation}\label{beta}
(i_{f}\IPS{(s,t)})f=(i_{f}(s^\top\wedge t))f=(f^\top \cdot (s^\top
t) )f=(sf)^\top tf.
\end{equation}
From this description, we see immediately that
\begin{equation}\label{invmuk}
\IPS{(x\cdot s,t)} = \IPS{(s,x^\top\cdot t)},\ \ \forall x\in
CL(V\oplus V^*).
\end{equation}
In particular $\IPS{(g\cdot s,g\cdot t)}=\pm\IPS{(s,t)}$ for any $g
\in \mathrm{Spin}(V\oplus V^*)$, yielding the following result.
\begin{prop}
The Mukai pairing is invariant under the identity component of Spin:
\begin{equation*}
\IPS{(g\cdot s,g\cdot t)}=\IPS{(s,t)}\ \ \forall\  g \in
\mathrm{Spin}_0(V\oplus V^*).
\end{equation*}
Therefore it determines a $Spin_0$-invariant bilinear form on
$S=\wedge^\bullet V^*\otimes(\det V)^{1/2}$.
\end{prop}
For example, we have $\IPS{(\exp B \cdot s,\exp B\cdot
s)}=\IPS{(s,t)}$, for any $B\in \wedge^2 V^*$.  The pairing is
non-degenerate, and is symmetric or skew-symmetric depending on the
dimension of $V$, since
\begin{equation*}
\IPS{(s,t)}=(-1)^{m(m-1)/2}\IPS{(t,s)}.
\end{equation*}
We also see from degree considerations that for $m$ even,
$\IPS{(S^+,S^-)}=0$.

\begin{example}
Suppose $V$ is 4-dimensional; then the Mukai pairing is symmetric,
and the even spinors are orthogonal to the odd spinors.  The inner
product of even spinors $\rho=\rho_0+\rho_2+\rho_4$ and
$\sigma=\sigma_0+\sigma_2+\sigma_4$ is given by
\begin{align*}
\IPS{(\rho,\sigma)}&=[(\rho_0 - \rho_2 + \rho_4)\wedge(\sigma_0 +
\sigma_2+ \sigma_4)]_4\\
&=\rho_0\wedge\sigma_4 - \rho_2\wedge\sigma_2 +
\rho_4\wedge\sigma_0.
\end{align*}
The inner product of odd spinors $\rho=\rho_1+\rho_3$ and
$\sigma=\sigma_1+\sigma_3$ is given by
\begin{align*}
\IPS{(\rho,\sigma)}&=[(\rho_1 - \rho_3)\wedge(\sigma_1 +
\sigma_3)]_4\\
&=\rho_1\wedge\sigma_3 - \rho_3\wedge\sigma_1.
\end{align*}
\end{example}

\subsection{Pure spinors and polarizations}\label{corres}

Let $\varphi\in\wedge^\bullet V^*$ be a nonzero spinor. We define
its \emph{null space} $L_\varphi<V\oplus V^*$, as follows:
\begin{equation*}
L_\varphi=\{v\in V\oplus V^*\ :\ v\cdot\varphi=0\},
\end{equation*}
and it is clear from~\eqref{repn} that $L_\varphi$ depends
equivariantly on $\varphi$, i.e.:
\begin{equation*}
L_{g\cdot\varphi}=\rho(g)L_\varphi\ \ \forall\ g\in
\mathrm{Spin}(V\oplus V^*).
\end{equation*}
The key property of null spaces is that they are isotropic: if
$v,w\in L_\varphi$, we have
\begin{equation*}
\IP{v,w}\varphi = \tfrac{1}{2}(vw+wv)\cdot\varphi=0.
\end{equation*}
\begin{defn}
A spinor $\varphi$ is called \emph{pure} when $L_\varphi$ is maximal
isotropic.
\end{defn}

Every maximal isotropic subspace $L\subset V\oplus V^*$ is
represented by a unique line $K_L\subset \wedge^\bullet V^*$ of pure
spinors, as we now describe. As we saw in section \ref{maxiso}, any
maximal isotropic $L(\Delta,\eps)$ may be expressed as the
$B$-transform of $L(\Delta,0)$ for $B$ chosen such that
$i^*B=-\veps$. The pure spinor line with null space
$L(\Delta,0)=\Delta + \Ann(\Delta)$ is precisely $\det (\Ann
(\Delta))\subset \wedge^k V^*$, for $k$ the codimension of
$\Delta\subset V$. Hence we obtain the following result.

\begin{prop}[\cite{Chevalley}, III.1.9.]\label{evenodd}
Let $L(\Delta,\veps)\subset V\oplus V^*$ be maximal isotropic, let
$(\theta_1,\ldots,\theta_k)$ be a basis for $\Ann(\Delta)$, and let
$B\in\wedge^2 V^*$ be a 2-form such that $i^*B=-\veps$, where
$i:\Delta\hookrightarrow V$ is the inclusion. Then the pure spinor
line $K_L$ representing $L(\Delta,\eps)$ is generated by
\begin{equation}\label{represent}
\varphi = \exp(B)\theta_1\wedge\cdots\wedge\theta_k.
\end{equation}
Note that $\varphi$ is of even or odd degree according as $L$ is of
even or odd parity.
\end{prop}
$L$ may also be described as a generalized graph on $V^*$
via~\eqref{othproj}, expressing it as a $\beta$-transform of
$L(F,0)$, which has associated pure spinor line $\det (L\cap V)$. As
a result we obtain the following complement to
Proposition~\ref{evenodd}.
\begin{prop}\label{bothways}
Given a subspace inclusion $i:\Delta\hookrightarrow V$ and a 2-form
$B\in \wedge^2 V^*$, there exists a bivector $\beta\in \wedge^2 V$,
such that
\[
e^B\det(\Ann (\Delta)) = e^\beta\det(\Ann(L\cap V))
\]
is an equality of pure spinor lines, where $L = L(\Delta,-i^*B)$.
Note that the image of $\beta$ in $\wedge^2(V/(L\cap V))$ is unique.
\end{prop}

The pure spinor line $K_L$ determined by $L$ forms the beginning of
an induced filtration on the spinors, as we now describe. Recall
that the Clifford algebra is a $\ZZ/2\ZZ$-graded, $\ZZ$-filtered
algebra with
\begin{equation*}
CL(V\oplus V^*)=CL^{2m}\supset CL^{2m-1}\supset\cdots\supset CL^0=\RR\\
\end{equation*}
where $CL^{k}$ is spanned by products of $\leq k$ generators. By the
Clifford action on $K_L$, we obtain a filtration of the spinors:
\begin{equation}\label{filtration}
K_L=F_0\subset F_1\subset \cdots\subset F_{m}=S\\
\end{equation}
where $F_k$ is defined as $CL^k\cdot K_L$. Note that $CL^k\cdot
K_L=CL^m\cdot K_L$ for $k>m$ since $L$ annihilates $K_L$.  Also,
using the inner product, we have the canonical isomorphism
$L^*=(V\oplus V^*)/L$, and so $F_k/F_{k-1}$ is isomorphic to
$\wedge^kL^*\otimes K_L$.

The filtration just described becomes a grading when a maximal
isotropic $L'\subset V\oplus V^*$ complementary to $L$ is chosen.
Then we obtain a $\ZZ$-grading on $S=\wedge^\bullet V^*$ of the form
(for $\dim V$ even, i.e. $m=2n$)
\begin{equation}\label{altgrad}
S = U^{-n}\oplus \cdots\oplus U^n,
\end{equation}
where $U^{-n}=K_L$ and $U^k = (\wedge^{k+n}L')\cdot K_L$, using the
inclusion as a subalgebra $\wedge^\bullet L'\subset CL(V\oplus
V^*)$.  Using the inner product to identify $L'=L^*$, we obtain the
natural isomorphism $U^k = \wedge^{k+n} L^*\otimes K_L$ and, summing
over $k$,
\begin{equation}\label{spineq}
\wedge^\bullet V^* = \wedge^\bullet L^*\otimes K_L.
\end{equation}
The line $U^n = \det L'\cdot U^{-n}$ is the pure spinor line
determining $L'$ and since $L\cap L'=\{0\}$, the Mukai pairing gives
an isomorphism
\begin{equation}\label{mukpai}
U^{-n}\otimes U^n= \det V^*,
\end{equation}
as can be seen from the nondegenerate pairing between $\wedge^0 V^*$
and $\det V^*$ and the Spin-invariance of the Mukai pairing.  This
is an example of how the Mukai pairing determines the intersection
properties of maximal isotropics, and it can be phrased as follows:
\begin{prop}[\cite{Chevalley}, III.2.4.]\label{chevalley}
Maximal isotropics $L,L'$ satisfy $L\cap L'=\{0\}$ if and only if
their pure spinor representatives $\varphi,\varphi'$ satisfy
$\IPS{(\varphi,\varphi')}\neq 0$.
\end{prop}
Furthermore, applying~\eqref{invmuk}, we see that the Mukai pairing
provides a nondegenerate pairing for all $k$:
\begin{equation*}
U^{-k}\otimes U^{k}\lra \det V^*
\end{equation*}

Finally, rewriting~\eqref{mukpai}, we have the isomorphism
$K_L\otimes \det L^*\otimes K_L = \det V^*$, which, combined
with~\eqref{spineq}, yields the canonical isomorphism
\begin{equation*}
\wedge^\bullet V^*\otimes(\det V)^{1/2}=\wedge^\bullet L^*\otimes
(\det L)^{1/2},
\end{equation*}
showing that tensoring with the half densities renders $S$
independent of polarization.

\subsection{The spin bundle for $T\oplus T^*$}

Consider the direct sum of the tangent and cotangent bundles
$T\oplus T^*$ of a real $m$-dimensional smooth manifold $M$. This
bundle is endowed with the same canonical bilinear form and
orientation we described on $V\oplus V^*$. Therefore, while $T\oplus
T^*$ is associated to a $GL(m)$ principal bundle, we may view it as
having natural structure group $SO(m,m)$.

It is well-known that an oriented bundle with Euclidean structure
group $SO(n)$ admits spin structure if and only if the second
Stiefel-Whitney class vanishes, i.e. $w_2(E)=0$.  For oriented
bundles with metrics of indefinite signature, we find the
appropriate generalization in~\cite{Karoubi}, which we now
summarize.

If an orientable bundle $E$ has structure group $SO(p,q)$, we can
always reduce the structure group to its maximal compact subgroup
$\mathrm{S}(\mathrm{O}(p)\times \mathrm{O}(q))$.  This reduction is
equivalent to the choice of a maximal positive definite subbundle
$C^+<E$, which allows us to write $E$ as the direct sum $E=C^+\oplus
C^-$, where $C^-=(C^+)^\perp$ is negative definite.
\begin{prop}[\cite{Karoubi}, 1.1.26]
The $SO(p,q)$ bundle $E$ admits $\mathrm{Spin}(p,q)$ structure if
and only if $w_2(C^+)=w_2(C^-)$.
\end{prop}
In the special case of $T\oplus T^*$, which has signature $(m,m)$,
the positive and negative definite bundles $C^\pm$ project
isomorphically via $\pi_T:T\oplus T^*\rightarrow T$ onto the tangent
bundle.  Hence the condition $w_2(C^+)=w_2(C^-)$ is always satisfied
for $T\oplus T^*$, yielding the following result.
\begin{prop}
The $SO(m,m)$ bundle $T\oplus T^*$ always admits a
$\mathrm{Spin}(m,m)$ structure.
\end{prop}

By the action defined in~\eqref{cliff}, the differential forms
$\wedge^\bullet T^*$ are a Clifford module for $T\oplus T^*$;
indeed, for an orientable manifold, we see immediately from the
decomposition~\eqref{spinors} of the spin representation under
$GL^+(m)$ that there is a natural choice of spin bundle, namely
\begin{equation*}
S = \wedge^\bullet T^*\otimes(\det T)^{1/2}.
\end{equation*}
The Mukai pairing may then either be viewed as a nondegenerate
pairing
\begin{equation*}
\wedge^\bullet T^*\otimes \wedge^\bullet T^*\lra \det T^*,
\end{equation*}
or as a bilinear form on the spinors $S$. In the rest of the paper
we will make frequent use of the correspondence between maximal
isotropics in $T\oplus T^*$ and pure spinor lines in $\wedge^\bullet
T^*$ to describe structures on $T\oplus T^*$ in terms of
differential forms.
\section{The Courant bracket}\label{courant}

The Courant bracket, introduced in~\cite{Courant},\cite{CourWein},
is an extension of the Lie bracket of vector fields to smooth
sections of $T\oplus T^*$. It differs from the Lie bracket in
certain important respects.  Firstly, its skew-symmetrization does
not satisfy the Jacobi identity; as shown in~\cite{RW} it defines
instead a Lie algebra up to homotopy. Secondly, it has an extended
symmetry group; besides the diffeomorphism symmetry it shares with
the Lie bracket, it admits B-field gauge transformations.  Finally,
as shown in~\cite{SeveraWeinstein}, the bracket may be `twisted' by
a closed 3-form $H$, and so may be viewed as naturally associated to
an $S^1$-gerbe.

Like the Lie bracket, the Courant bracket may be defined as a
derived bracket; we describe this following the treatment
in~\cite{KS}. We also present, following~\cite{Roytenberg}, certain
characteristic properties of the Courant bracket as formalized
by~\cite{LWX} in the notion of Courant algebroid, a generalization
of Lie algebroid.

\subsection{Derived brackets}

The interior product of vector fields $X\mapsto i_X$ defines an
effective action of vector fields on differential forms by
derivations of degree $-1$.  Taking supercommutator with the
exterior derivative, we obtain a derivation of degree $0$, the Lie
derivative:
\begin{equation*}
\LL_X = [d,i_X] = di_X + i_Xd.
\end{equation*}
Taking commutator with another interior product, we obtain an
expression for the Lie bracket:
\begin{equation*}
i_{[X,Y]} = [\LL_X,i_Y].
\end{equation*}
In this sense, the Lie bracket on vector fields is `derived' from
the Lie algebra of endomorphisms of $\Omega^\bullet(M)$.  Using the
spinorial action~\eqref{cliff} of $T\oplus T^*$ on forms, we may
define the Courant bracket of sections $e_i\in C^\infty(T\oplus
T^*)$ in the same way:
\begin{equation}\label{derivbrack}
[e_1,e_2]\cdot \varphi = [[d,e_1\cdot],e_2\cdot]\varphi\ \ \ \
\forall \varphi\in \Omega^\bullet(M).
\end{equation}
Although this bracket is not skew-symmetric, it follows from the
fact that $d^2=0$ that the following Jacobi identity holds:
\begin{equation}\label{jacobi}
[[e_1,e_2],e_3]=[e_1,[e_2,e_3]] - [e_2,[e_1,e_3]].
\end{equation}

As observed in~\cite{SeveraWeinstein}, one may replace the exterior
derivative in~\eqref{derivbrack} by the twisted operator $d_H\varphi
= d\varphi + H\wedge \varphi$ for a real 3-form $H\in \Omega^3(M)$.
The resulting derived bracket then satisfies
\begin{equation*}
[[e_1,e_2],e_3]=[e_1,[e_2,e_3]] - [e_2,[e_1,e_3]] +
i_{\pi{e_3}}i_{\pi{e_2}}i_{\pi{e_1}}dH,
\end{equation*}
where $\pi:T\oplus T^*\lra T$ is the first projection.  When the
last term vanishes, the bracket is called a Courant bracket:
\begin{defn}
Let $e_1,e_2\in C^\infty(T\oplus T^*)$. Then their $H$-twisted
Courant bracket $[e_1,e_2]_H\in C^\infty(T\oplus T^*)$ is defined by
the expression
\begin{equation*}
[e_1,e_2]_H\cdot \varphi = [[d_H,e_1\cdot],e_2\cdot]\varphi\ \ \ \
\forall \varphi\in \Omega^\bullet(M),
\end{equation*}
where $d_H\varphi = d\varphi + H\wedge\varphi$ for $H$ a real closed
3-form.
\end{defn}
Expanding this expression for $e_1=X+\xi$ and $e_2=Y+\eta$, we
obtain
\begin{equation}\label{stdcour}
[X+\xi,Y+\eta]_H = [X,Y]  +  \LL_X\eta - i_Yd\xi + i_Xi_YH.
\end{equation}
The Courant bracket, for any closed 3-form $H$, satisfies the
following conditions~\cite{Roytenberg}:
\begin{itemize}
\item[C1)] $\Cour{e_1,\Cour{e_2,e_3}} = \Cour{\Cour{e_1,e_2},e_3}
+ \Cour{e_2,\Cour{e_1,e_3}}$,
\item[C2)] $\pi(\Cour{e_1,e_2})=[\pi(e_1),\pi(e_2)]$,
\item[C3)] $\Cour{e_1,fe_2}=f\Cour{e_1,e_2}+ (\pi(e_1)f) e_2,\ f\in
C^\infty(M)$,
\item[C4)] ${\pi(e_1)}\IP{e_2,e_3}= \IP{\Cour{e_1,e_2},e_3}
+ \IP{e_2, \Cour{e_1, e_3}}$,
\item[C5)] $[e_1,e_1] = \pi^*d\IP{e_1,e_1}$.
\end{itemize} In~\cite{LWX}, these properties were promoted to axioms defining
the notion of \emph{Courant algebroid}, which is a real vector
bundle $E$ equipped with a bracket $[\cdot,\cdot]$, nondegenerate
inner product $\IP{\cdot,\cdot}$, and bundle map $\pi:E\lra T$
(called the \emph{anchor}) satisfying the conditions C1)--C5) above.

Note that if the bracket were skew-symmetric, then axioms C1)--C3)
would define the notion of \emph{Lie algebroid}; axiom C5) indicates
that the failure to be a Lie algebroid is measured by the inner
product, which itself is invariant under the adjoint action by axiom
C4).

The Courant algebroid structure on $T\oplus T^*$ has surjective
anchor with isotropic kernel given by $T^*$; such Courant algebroids
are called \emph{exact}.
\begin{defn}
A Courant algebroid $E$ is \emph{exact} when the following sequence
is exact:
\begin{equation}\label{exact}
\xymatrix{0\ar[r]&T^*\ar[r]^{\pi^*}&E\ar[r]^\pi&T\ar[r]&0},
\end{equation}
where $E$ is identified with $E^*$ using the inner product.
\end{defn}
Exactness at the middle place forces $\pi^*(T^*)$ to be isotropic
and therefore the inner product on $E$ must be of split signature.
It is then always possible to choose an isotropic splitting $s:T\lra
E$ for $\pi$, yielding an isomorphism $E\cong T\oplus T^*$ taking
the Courant bracket to that given by~\eqref{stdcour}, where $H$ is
the curvature of the splitting, i.e.
\begin{equation}\label{splitcurv}
i_Xi_Y H = s^*\Cour{s(X),s(Y)}, \ \ \ X,Y\in T.
\end{equation}
Isotropic splittings of $\eqref{exact}$ are acted on transitively by
the 2-forms $B\in \Omega^2(M)$ via transformations of the form
$X+\xi\mapsto X+\xi + i_XB$, or more invariantly,
\begin{equation*}
e\mapsto e + \pi^*i_{\pi e}B,\ \ \ e\in E.
\end{equation*}
Such a change of splitting modifies the curvature $H$ by the exact
form $dB$. Hence the cohomology class $[H]\in H^3(M,\mathbb{R})$,
called the \textit{\v{S}evera class}, is independent of the
splitting and determines the exact Courant algebroid structure on
$E$ completely.

\subsection{Symmetries of the Courant bracket}

The Lie bracket of smooth vector fields is invariant under
diffeomorphisms; in fact, there are no other symmetries of the
tangent bundle preserving the Lie bracket.
\begin{prop}
Let $F$ be an automorphism of the tangent bundle covering the
diffeomorphism $\varphi:M\lra M$.  If $F$ preserves the Lie bracket,
then $F=\varphi_*$, the derivative of $\varphi$.
\end{prop}
\begin{proof}
Let $G=\varphi_*^{-1}\circ F$, so that it is an automorphism of the
Lie bracket covering the identity map.  Then, for any vector fields
$X,Y$ and $f\in C^\infty(M)$ we have $G([fX,Y])=[G(fX),G(Y)]$.
Expanding, we see that $Y(f)G(X)=G(Y)(f)G(X)$ for all $X,Y,f$. This
can only hold when $G(Y)=Y$ for all vector fields $Y$, i.e. $G=Id$,
yielding finally that $F=\varphi_*$.
\end{proof}

Since the Courant bracket on $T\oplus T^*$ depends on a 3-form $H$,
it may appear at first glance to have a smaller symmetry group than
the Lie bracket. However, as was observed in~\cite{SeveraWeinstein},
the spinorial action of 2-forms~\eqref{Bfield} satisfies
\begin{equation*}
e^{-B}d_H e^B = d_{H+dB},
\end{equation*}
and therefore we obtain the following action on derived brackets:
\begin{equation}\label{nonclosed}
e^B[e^{-B}\cdot,e^{-B}\cdot]_H = [\cdot,\cdot]_{H+dB}.
\end{equation}
We see immediately from~\eqref{nonclosed} that closed 2-forms act as
symmetries of any exact Courant algebroid.
\begin{defn}
A B-field transformation is the automorphism of an exact Courant
algebroid $E$ defined by a closed 2-form $B$ via
\begin{equation*}
e\mapsto e + \pi^*i_{\pi e}B.
\end{equation*}
\end{defn}
A diffeomorphism $\varphi:M\lra M$ lifts to an orthogonal
automorphism of $T\oplus T^*$ given by
\begin{equation*}
\left(\begin{matrix}\varphi_*&0\\0&{\varphi^*}^{-1}\end{matrix}\right),
\end{equation*}
which we will denote by $\varphi_*$.  It acts on the Courant bracket
via
\begin{equation}\label{diffact}
\varphi_*[\varphi_*^{-1}\cdot,\varphi_*^{-1}\cdot]_H =
[\cdot,\cdot]_{{\varphi^*}^{-1}H}.
\end{equation}
Combining~\eqref{diffact} with~\eqref{nonclosed}, we see that the
composition $F = \varphi_*e^B$ is a symmetry of the $H$-twisted
Courant bracket if and only if ${\varphi^*}H -H= dB$.  We now show
that such symmetries exhaust the automorphism group.
\begin{prop}\label{formiso}
Let $F$ be an orthogonal automorphism of $T\oplus T^*$, covering the
diffeomorphism $\varphi:M\lra M$, and preserving the $H$-twisted
Courant bracket. Then $F = \varphi_*e^B$ for a unique 2-form
$B\in\Omega^2(M)$ satisfying $\varphi^*H - H=dB$.
\end{prop}
\begin{proof}
Let $G=\varphi_*^{-1}\circ F$, so that it is an automorphism of
$T\oplus T^*$ covering the identity satisfying
$G[G^{-1}\cdot,G^{-1}\cdot]_H = [\cdot,\cdot]_{\varphi^*H}$.  In
particular, for any sections $x,y\in C^\infty(T\oplus T^*)$ and
$f\in C^\infty(M)$ we have $G[x,fy]_H=[Gx,Gfy]_{\varphi^*H}$, which,
using axiom $C3)$, implies
\begin{align*}
\pi(x)(f) Gy = \pi(Gx)(f) Gy.
\end{align*}
Therefore, $\pi G = \pi$, and so $G$ is an orthogonal map preserving
projections to $T$. This forces it to be of the form $G=e^B$, for
$B$ a uniquely determined 2-form. By~\eqref{nonclosed}, $B$ must
satisfy $\varphi^*H-H=dB$. Hence we have $F = \varphi_*e^B$, as
required.
\end{proof}

An immediate corollary of this result is that the automorphism group
of an exact Courant algebroid $E$ is an extension of the
diffeomorphisms preserving the cohomology class $[H]$ by the abelian
group of closed 2-forms:
\begin{equation*}\label{extenex}
\xymatrix{0\ar[r]&\Omega^2_{cl}(M)\ar[r]&\mathrm{Aut}(E)\ar[r]&\Diff_{[H]}(M)\ar[r]&0}.
\end{equation*}

Derivations of a Courant algebroid $E$ are linear first order
differential operators $D_X$ on $C^\infty(E)$, covering vector
fields $X$ and satisfying
\begin{equation*}
\begin{split}
X\IP{\cdot,\cdot} = \IP{D_X\ \cdot,\cdot}+\IP{\cdot,D_X\ \cdot},\\
D_X[\cdot,\cdot]=[D_X\ \cdot,\cdot]+[\cdot,D_X\ \cdot].
\end{split}
\end{equation*}
Differentiating a 1-parameter family of automorphisms $F_t =
\varphi_*^t e^{B_t}$, $F_0=\mathrm{Id}$, and using the convention
for Lie derivative
\[
\LL_X = -\tfrac{d}{dt}\big|_{t=0}\varphi_*^t,
\]
we see that the Lie algebra of derivations of the $H$-twisted
Courant bracket consists of pairs $(X,b)\in C^\infty(T)\oplus
\Omega^2(M)$ such that $\LL_X H= db$, which act via
\begin{equation}\label{deract}
(X,b)\cdot(Y+\eta) = \LL_X(Y+\eta) - i_Yb.
\end{equation}
Therefore the algebra of derivations of an exact Courant algebroid
algebroid $E$ is an abelian extension of the smooth vector fields by
the closed 2-forms:
\begin{equation*}
\xymatrix{0\ar[r]&\Omega^2_{cl}(M)\ar[r]&\mathrm{Der}(E)\ar[r]&C^\infty(T)\ar[r]&0}.
\end{equation*}
It is clear from axioms $\mathrm{C1), C4)}$ that the left adjoint
action $\ad_v:w\mapsto[v,w]$ defines a derivation of the Courant
algebroid. However, $\ad$ is neither surjective nor injective;
rather, for $E$ exact, it induces the following exact sequence:
\begin{equation}\label{intder}
\xymatrix{0\ar[r]&\Omega^1_{cl}(M)\ar[r]^{\pi^*}&C^\infty(E)\ar[r]^{\ad}&\mathrm{Der}(E)\ar[r]^\chi&H^2(M,\RR)\ar[r]&0},
\end{equation}
where $\Omega^1_{cl}(M)$ denotes the closed 1-forms and we define
$\chi(D_X)=[i_XH-b]\in H^2(M,\RR)$ for $D_X = (X,b)$ as above.

The image of $\ad$ in sequence~\eqref{intder} defines a Lie
subalgebra of $\mathrm{Der}(E)$, and so suggests the definition of a
subgroup of the automorphism group of $E$ analogous to the subgroup
of Hamiltonian symplectomorphisms.
\begin{prop}
Let $D_{X_t}=(X_t,b_t)\in C^\infty(T)\oplus \Omega^2(M)$ be a
(possibly time-dependent) derivation of the $H$-twisted Courant
bracket on a compact manifold, so that it satisfies
$\LL_{X_t}H=db_t$ and acts via~\eqref{deract}. Then it generates a
1-parameter subgroup of Courant automorphisms
\begin{equation}\label{oneflow}
F^t_{D_X} = \varphi^t_*e^{B_t},\ \ t\in\RR,
\end{equation}
where $\varphi^t$ denotes the flow of the vector field $X_t$ for a
time $t$ and
\begin{equation}\label{intbee}
B_t = \int_{0}^{t} \varphi^*_sb_s \ ds.
\end{equation}
\end{prop}
\begin{proof}
First we see that $F^t_{D_X}$ is indeed an automorphism, since
\begin{align*}
dB_t &= \int_{0}^t\varphi_s^*(\LL_{X_s}H)\ ds=\tfrac{d}{du}\big|_{u=0}\int_{0}^t\varphi_s^*\varphi_{u}^*H\ ds\\
&=\tfrac{d}{du}\big|_{u=0}\int_{u}^{t+u}\varphi_{s'}^*H\ ds'\\
&= \varphi^*_{t}H - H,
\end{align*}
which proves the result by Proposition~\ref{formiso}.  To see that
it is a 1-parameter subgroup, observe that $e^B\varphi_* = \varphi_*
e^{\varphi^*B}$ for any $\varphi\in\Diff(M)$ and $B\in\Omega^2(M)$,
so that
\begin{align*}
\varphi^{t_1}_*e^{B_{t_1}}\varphi^{t_2}_*e^{B_{t_2}}=
\varphi_*^{t_1+t_2}e^{\varphi_{t_2}^*B_{t_1} +
B_{t_2}}=\varphi_*^{t_1+t_2}e^{B_{t_1+t_2}},
\end{align*}
where we use the expression~\eqref{intbee} for the final equality.
\end{proof}
Certain derivations $(X,b)$ are in the kernel of $\chi$
in~\eqref{intder}, namely those for which $b = i_XH+ d\xi$ for a
1-form $\xi$; we call these \emph{exact} derivations.  A smooth
1-parameter family of automorphisms  $F_t = \varphi_*^te^{B_t}$ from
$F_0=\id$ to $F_1= F$ is called an \emph{exact} isotopy when it is
generated by a smooth time-dependent family of exact derivations.
\begin{defn}\label{exaut}
An automorphism $F\in\mathrm{Aut}(E)$ is called \emph{exact} if
there is an exact isotopy $F_t$ from $F_0=\id$ to $F_1=F$.  This
defines the subgroup of exact automorphisms of any Courant
algebroid:
\[
\mathrm{Aut}_{ex}(E)\subset \mathrm{Aut}(E).
\]
\end{defn}
If $\ad(v_t)$ generates the exact isotopy $F_t$ and $F$ is any
automorphism, then the conjugation $FF_t F^{-1}$ is also an exact
isotopy, generated by the family of derivations
\[
\ad(F(v_t)).
\]
Therefore we see that $\mathrm{Aut}_{ex}(E)$ is a \emph{normal}
subgroup, in analogy with the group of Hamiltonian
symplectomorphisms.

\subsection{Relation to $S^1$-gerbes}\label{gerb}

The Courant bracket is part of a hierarchy of brackets on the
bundles $T\oplus \wedge^p T^*$, $p=0,1,\ldots$, defined by the same
formula
\begin{equation*}
[X+\sigma,Y+\tau]=[X,Y]+\LL_X\tau-i_Yd\sigma + i_Xi_Y F,
\end{equation*}
where now $\sigma,\tau\in C^\infty(\wedge^p T^*)$ and $F$ is a
closed form of degree $p+2$.

For $p=0$, the Courant bracket on $T\oplus 1$ is given by
\begin{equation}\label{pequals0}
[X+f,Y+g]=[X,Y]+Xg-Yf + i_Xi_YF.
\end{equation}
When $F/2\pi$ is integral, we recognize this bracket as coming from
the Atiyah Lie algebroid $\mathcal{A}=TP/S^1$ associated to a
principal $S^1$ bundle $\pi:P\lra M$, i.e.
\begin{equation}\label{atiy}
\xymatrix{0\ar[r]&1\ar[r]&\mathcal{A}\ar[r]^{\pi_*}&TM\ar[r]& 0},
\end{equation}
where $1=\wedge^0T^*M$ denotes the trivial line bundle over $M$.  A
splitting of this sequence provides a connection 1-form
$\theta\in\Omega^1(P)$ with curvature $d\theta = \pi^*F$, and we see
that the natural Lie bracket on $S^1$-invariant vector fields may be
written, in horizontal and vertical components, as
\begin{equation*}
[X_h+f\partial_t,Y_h+g\partial_t]=[X,Y]_h + (Xg-Yf +
i_Xi_YF)\partial_t,
\end{equation*}
where $X\mapsto X_h$ denotes horizontal lift and $\partial_t$ is the
vector field generating the principal $S^1$ action.  In this way we
recover the expression~\eqref{pequals0}.  The symmetries of
$\mathcal{A}$ covering the identity consist of closed 1-forms $A$
acting via $X+f\mapsto X+f+i_XA$, which may be interpreted as the
action of tensoring with a trivialization of a flat unitary line
bundle. When $[A]\in H^1(M,\ZZ)$, then it represents the action of a
gauge transformation on $P$, modulo constant gauge transformations.

Just as a Lie algebroid of the form~\eqref{atiy} may be associated
with a $S^1$ bundle when $[F]/2\pi$ is integral, an exact Courant
algebroid may be associated with a $S^1$ gerbe when $[H]/2\pi$ is
integral.  A $S^1$ gerbe may be specified, given an open cover
$\{U_i\}$, by complex line bundles $L_{ij}$ on $U_{i}\cap U_j$ with
isomorphisms $L_{ji}\cong L_{ij}^*$, and trivializations
$\theta_{ijk}$ of $L_{ij}\otimes L_{jk}\otimes L_{ki}$ such that
$\delta\theta = 1$ in the \v{C}ech complex.  A
0-connection~\cite{Chatterjee} on a gerbe is then specified by
choosing connections $\nabla_{ij}$ on $L_{ij}$ such that the induced
connection on threefold intersections obeys $\nabla\theta = 0$.
Letting $F_{ij}$ be the curvature 2-forms of $\nabla_{ij}$, this
implies that $\delta F = 0$, i.e.:
\begin{equation}\label{coc}
F_{ij} + F_{jk} + F_{ki} = 0.
\end{equation}
As explained in~\cite{Hitchin}, we may then construct an exact
Courant algebroid $E$ as follows: glue $T\oplus T^*$ over $U_i$ with
$T\oplus T^*$ over $U_j$ using the transition function $\Phi_{ij}$
given by the $B$-field transform
\begin{equation*}
\Phi_{ij}= e^{F_{ij}}.
\end{equation*}
By the cocycle condition~\eqref{coc}, we see that $\delta \Phi = 1$
and therefore it defines a bundle $E$.  Since $\Phi_{ij}$ preserves
projections to $T$, we see that $E$ is an extension of $T$ by $T^*$,
as required. Finally, equipping $E|_{U_i}$ with the standard Courant
bracket ($H=0$) and inner product, we observe that since $\Phi_{ij}$
is orthogonal and preserves the Courant bracket, $E$ inherits the
structure of an exact Courant algebroid.  Choosing a splitting of
the exact sequence then corresponds to the choice of $1$-connection
for the gerbe, i.e. the choice of local 2-forms $B_i$ such that $F =
\delta B$.  This then determines the global curvature 3-form of the
gerbe $H = dB_i$, for which $H/2\pi$ is integral.  In this sense, an
exact Courant algebroid with integral curvature may be viewed as the
generalized Atiyah sequence of a $S^1$ gerbe.

Similarly to the case of $p=0$, the symmetries of $E$ covering the
identity consist of closed $2$-forms $B$ acting via $B$-field
transforms, and these may be interpreted as trivializations of a
flat gerbe. The difference of two such trivializations, a line
bundle with connection, acts as a gauge transformation (integral
$B$-field).

In the case that $F/2\pi$ is not integral, the Lie
algebroid~\eqref{atiy} may be interpreted as the Atiyah sequence of
a trivialization of a $S^1$ gerbe with flat connection; similarly, a
general exact Courant algebroid may be associated with a
trivialization of a $S^1$ 2-gerbe with flat connection.  The fact
that such trivializations may be tensored together accounts for the
Baer sum operation described by \v{S}evera (see
\cite{Sevlet},\cite{Bressler}) on exact Courant algebroids.

\subsection{Dirac structures}

The Courant bracket fails to be a Lie bracket due to exact terms
involving the inner product $\IP{,}$.  Therefore, upon restriction
to a subbundle $L\subset T\oplus T^*$ which is involutive (closed
under the Courant bracket) as well as being isotropic, the anomalous
terms vanish.  Then $(L,[,],\pi)$ defines a Lie algebroid, with
associated differential graded algebra $(\wedge^\bullet L^*, d_L)$,
just as the de Rham complex is associated to the canonical Lie
algebroid structure on the tangent bundle.

In fact, the Courant bracket itself places a tight constraint on
which proper subbundles may be involutive a priori:
\begin{prop}\label{isot}
If $L\subset E$ is an involutive subbundle of an exact Courant
algebroid, then $L$ must be isotropic, or of the form
$\pi^{-1}(\Delta)$, for $\Delta$ an integrable distribution in $T$.
\end{prop}
\begin{proof}
Suppose that $L\subset E$ is involutive, but not isotropic, i.e.
there exists $v\in C^\infty(L)$ such that $\IP{v,v}\neq 0$ at some
point $m\in M$.  Then for any $f\in C^\infty(M)$,
\begin{align*}
[fv,v] = f[v,v] - (\pi(v)f)v + 2\IP{v,v}df,
\end{align*}
implying that $df|_m\in L|_m$ for all $f$, i.e. $T^*|_m\subset
L|_m$. Since $T^*|_m$ is isotropic, this inclusion must be proper,
i.e. $L|_m=\pi^{-1}(\Delta|_m)$, where $\Delta=\pi(L)$ is nontrivial
at $m$.  Hence the rank of $L$ must exceed the maximal dimension of
an isotropic subbundle. This implies that $T^*|_m<L|_m$ at every
point $m$, and hence that $\Delta$ is a smooth subbundle of $T$,
which must itself be involutive. Hence $L=\pi^{-1}(\Delta)$, as
required.
\end{proof}

\begin{defn}[Dirac structure]
A maximal isotropic subbundle $L\subset E$ of an exact Courant
algebroid is called an almost \emph{Dirac} structure.  If $L$ is
involutive, then the almost Dirac structure is said to be
integrable, or simply a \emph{Dirac} structure.
\end{defn}
The integrability of an almost Dirac structure $L$ may be expressed,
following \cite{Courant}, as the vanishing of the operator
$T_L(e_1,e_2,e_3)=\IP{[e_1,e_2],e_3}$ on sections of $L$. Using the
Clifford action we see that, for $\varphi$ a local generator of the
pure spinor line $K_L\subset \wedge^*T^*$ representing $L$, and
sections $e_i\in C^\infty(L)$,
\begin{equation}\begin{split}\label{tcour}
\IP{[e_1,e_2]_H,e_3}\varphi&=[[[d_H,e_1],e_2],e_3]\varphi\\
&=e_1\cdot e_2\cdot e_3\cdot d_H\varphi.
\end{split}\end{equation}
Therefore, as shown by Courant, integrability of a Dirac structure
is determined by the vanishing of a tensor $T_L\in C^\infty(\wedge^3
L^*)$.

We see from Proposition~\ref{Eeps} that at a point $p$, a Dirac
structure $L\subset T\oplus T^*$ has a unique description as a
generalized graph $L(\Delta,\veps)$, where $\Delta=\pi(L)$ is the
projection to $T$ and $\veps\in \wedge^2 \Delta^*$. Assuming that
$L$ is \emph{regular} near $p$ in the sense that $\Delta$ has
constant rank near $p$, we have the following description of the
integrability condition:
\begin{prop}\label{regular}
Let $\Delta\subset T$ be a subbundle and $\veps\in C^\infty(\wedge^2
\Delta^*)$.  Then the almost Dirac structure $L(\Delta,\veps)$ is
integrable for the $H$-twisted Courant bracket if and only if
$\Delta$ integrates to a foliation and $d_\Delta\veps=i^*H$, where
$d_\Delta$ is the leafwise exterior derivative.
\end{prop}
\begin{proof}
Let $i:\Delta\hookrightarrow T$ be the inclusion. Then
$d_\Delta:C^\infty(\wedge^k \Delta^*)\rightarrow
C^\infty(\wedge^{k+1}\Delta^*)$ is defined by $i^*\circ d =
d_\Delta\circ i^*$.
Suppose that
$X+\xi,Y+\eta\in C^\infty(L)$, i.e. $i^*\xi=i_X\veps$ and
$i^*\eta=i_Y\veps$. Consider the bracket $Z+\zeta=[X+\xi,Y+\eta]$;
if $L$ is Courant involutive, then $Z=[X,Y]\in C^\infty(\Delta)$,
showing $\Delta$ is involutive, and the difference
\begin{align*}
i^*\zeta-i_Z\veps &=i^*(\LL_X\eta - i_Yd\xi +
i_Xi_YH)-i_{[X,Y]}\veps\\
&= d_\Delta i_Xi_Y\veps + i_Xd_\Delta i_Y\veps - i_Yd_\Delta i_X\veps + i_Xi_Yi^*H -[[d_\Delta, i_X],i_Y]\veps\\
&=i_Yi_X(d_\Delta\veps - i^*H)
\end{align*}
must vanish for all $X+\xi,Y+\eta\in C^\infty(L)$, showing that
$d_\Delta\veps=i^*H$. Reversing the argument we see that the
converse holds as well.
\end{proof}
A consequence of this is that in a regular neighbourhood, a
$d_H$-closed generator for the canonical line bundle may always be
chosen:
\begin{corollary}\label{locclos}
Let $(\Delta,\veps)$ be as above and assume $L(\Delta,\veps)$ is
integrable; then for $B\in C^\infty(\wedge^2 T^*)$ such that $i^*B =
-\veps$, there exists a basis of sections
$(\theta_1,\ldots,\theta_k)$ for $\Ann(\Delta)$ such that
\[
\varphi = e^B\theta_1\wedge\cdots\wedge\theta_k
\]
is a $d_H$-closed generator for the pure spinor line $K_L$.
\end{corollary}
\begin{proof}
Let $\Omega = \theta_1\wedge\ldots\wedge\theta_k$.  By
Proposition~\ref{regular}, $\Delta$ is integrable, so
$(\theta_1,\ldots,\theta_k)$ can be chosen such that $d\Omega=0$.
Then we have
\[
d_H(e^B\Omega)= (dB + H)\wedge e^B\Omega=0,
\]
where the last equality holds since
$i^*(dB+H)=-d_\Delta\eps+i^*H=0$.
\end{proof}

In neighbourhoods where $\Delta$ is not regular, one may not find
$d_H$-closed generators for $K_L$; nevertheless, one has the
following useful description of the integrability condition.  Recall
that an almost Dirac structure $L\subset T\oplus T^*$ determines a
filtration \eqref{filtration} of the forms $\wedge^\bullet T^*$.
By~\eqref{tcour}, we see that $d_H$ takes $C^\infty(F_0)$ into
$C^\infty(F_3)$.  The integrability of $L$, however, requires that
$d_H$ take $C^\infty(F_0)$ into $C^\infty(F_1)$, as we now show.
\begin{theorem}\label{formint}
The almost Dirac structure $L\subset T\oplus T^*$ is involutive for
the $H$-twisted Courant bracket if and only if
\begin{equation}\label{dfilt}
d_H(C^\infty(F_0))\subset C^\infty(F_1),
\end{equation}
i.e. for any local trivialization $\varphi$ of $K_L$, there exists a
section $X+\xi\in C^\infty(T\oplus T^*)$ such that
\begin{equation*}
d_H\varphi = i_X\varphi+ \xi\wedge\varphi.
\end{equation*}
Furthermore, condition~\eqref{dfilt} implies that
\begin{equation}\label{dfiltk}
d_H(C^\infty(F_k))\subset C^\infty(F_{k+1})
\end{equation}
for all $k$.
\end{theorem}
\begin{proof}
Let $\varphi$ be a local generator for $K_L=F_0$. Then for
$e_1,e_2\in C^\infty(L)$, we have
\[\begin{split}
[e_1,e_2]_H\cdot \varphi &= [[d_H,e_1],e_2]\varphi \\
&=e_1\cdot e_2\cdot d_H\varphi,
\end{split}\]
and therefore $L$ is involutive if and only if $d_H\varphi$ is
annihilated by all products $e_1e_2$, $e_i\in C^\infty(L)$.  Since
$F_k$ is precisely the subbundle annihilated by products of $k+1$
sections of $L$, we obtain $d_H\varphi\in F_1$.  The subbundle $F_1$
decomposes in even and odd degree parts as $F_1 = F_0\oplus (T\oplus
T^*)\cdot F_0$, and since $d_H$ is of odd degree, we see that
$d_H\varphi\in (T\oplus T^*)\cdot K_L$, as required.  To
prove~\eqref{dfiltk}, we proceed by induction on $k$; let $\psi\in
F_k$, then since $[e_1,e_2]_H\cdot\psi = [[d_H,e_1],e_2]\psi$, we
have
\[
e_1\cdot e_2\cdot d_H\psi = d_H(e_1\cdot e_2\cdot \psi) + e_1\cdot
d_H (e_2\cdot\psi) - e_2\cdot d_H( e_1\cdot\psi) -
[e_1,e_2]_H\cdot\psi.
\]
All terms on the right hand side are in $F_{k-1}$ by induction,
implying that $d_H\psi\in F_{k+1}$, as required.
\end{proof}

Since the inner product provides a natural identification $(T\oplus
T^*)\cdot K_L = L^*\otimes K_L$, the previous result shows that the
pure spinor line generating a Dirac structure is equipped with an
operator
\begin{equation}\label{modla}
d_H: C^\infty(K_L)\lra C^\infty(L^*\otimes K_L),
\end{equation}
which satisfies $d_H^2=0$ upon extension to $C^\infty(\wedge^k
L^*\otimes K_L)$. This makes $K_L$ a Lie algebroid module for $L$,
i.e. a module over the differential graded Lie algebra
$(\wedge^\bullet L^*, d_L)$ associated to the Lie algebroid $L$
(see~\cite{Fernandes} for discussion of Lie algebroid modules).

\begin{example}
The cotangent bundle $T^*\subset T\oplus T^*$ is a Dirac structure
for any twist $H\in \Omega^3_{cl}(M)$.
\end{example}

\begin{example}\label{symp}
The tangent bundle $T\subset T\oplus T^*$ is itself maximal
isotropic and involutive for the Courant bracket with $H=0$, hence
defines a Dirac structure. Applying any 2-form $B\in \Omega^2(M)$,
we obtain
\begin{equation*}
\Gamma_B=e^B(T)=\{X+i_XB \ :\ X\in T\},
\end{equation*}
which is a Dirac structure for the $dB$-twisted Courant bracket.
Indeed, $T^*$ has no complementary Dirac structure unless $[H]=0$.
\end{example}
\begin{example}[Twisted Poisson geometry]\label{pois}
Applying a bivector field $\beta\in C^\infty(\wedge^2 T)$ as in
Example~\ref{betatrans} to the Dirac structure $T^*$, we obtain
\begin{equation}\label{poisson}
\Gamma_\beta=e^\beta(T^*)=\{i_\xi\beta+\xi \ :\ \xi\in T^*\}.
\end{equation}
As shown in~\cite{SeveraWeinstein}, this almost Dirac structure is
integrable with respect to the $H$-twisted Courant bracket if and
only if
\[
[\beta,\beta]=\wedge^3\beta^*(H),
\]
where the bracket denotes the Schouten bracket of bivector fields.
Such a structure is called a \emph{twisted Poisson structure}, and
becomes a usual Poisson structure when $\wedge^3\beta^*(H)=0$.
\end{example}
\begin{example}[Foliations]
Let $\Delta\subset T$ be a smooth distribution of constant rank.
Then the maximal isotropic subbundle
\begin{equation*}
\Delta\oplus \Ann(\Delta)\subset T\oplus T^*
\end{equation*}
is Courant involutive if and only if $\Delta$ is integrable and
$H|_\Delta = 0$.
\end{example}
\begin{example}[Complex geometry]\label{cxdirac}
An almost complex structure $J\in\End(T)$ determines a complex
distribution, given by the $-i$-eigenbundle $T_{0,1}<T\otimes\CC$ of
$J$.  Forming the maximal isotropic subbundle
\begin{equation*}
L_J=T_{0,1}\oplus\Ann(T_{0,1})=T_{0,1}\oplus T^*_{1,0},
\end{equation*}
we see from Proposition~\ref{regular} that
$L_J$ is integrable if and only if $T_{0,1}$ is involutive and
$i^*H=0$ for the inclusion $i:T^{0,1}\hookrightarrow T\otimes\CC$,
i.e. $H$ is of type $(1,2)+(2,1)$. Viewing $H$ as the curvature of a
gerbe, this means that the gerbe inherits a holomorphic structure
compatible with the underlying complex manifold. In this way,
integrable complex structures equipped with holomorphic gerbes can
be described by (complex) Dirac structures.
\end{example}

We may apply Theorem~\ref{formint} to give a simple description of
the \emph{modular vector field} of a Poisson structure (we
follow~\cite{LuEv}; for the case of twisted Poisson structures,
see~\cite{KoL}). The Dirac structure~\eqref{poisson} associated to a
Poisson structure $\beta$ has corresponding pure spinor line
generated by $\varphi=e^\beta\cdot v$, where $v\in C^\infty(\det
T^*)$ is a volume form on the manifold, which we assume to be
orientable.  By Theorem~\ref{formint}, there exists $X+\xi\in
C^\infty(T\oplus T^*)$ such that $d\varphi = (X+\xi)\cdot \varphi$.
Since $L_\beta$ annihilates $\varphi$, there is a unique $X_v\in
C^\infty(T)$, called the \emph{modular vector field} associated to
$(\beta,v)$, such that
\begin{equation}\label{mvf}
d\varphi = X_v\cdot\varphi.
\end{equation}
We see from applying $d$ to~\eqref{mvf} that
\[
\LL_{X_v}\varphi = d(X_v\cdot\varphi) + X_v\cdot X_v\cdot\varphi=0,
\]
implying immediately that $X_v$ is a Poisson vector field (i.e.
$[\beta,X_v]=0$) preserving the volume form $v$.  Of course the
modular vector field is not an invariant of the Poisson structure
alone; for $f\in C^\infty(M,\RR)$, one obtains
\[
X_{(e^f v)} = X_v + [\beta,f].
\]
As a result we see, following Weinstein~\cite{Weinmod}, that $X_v$
defines a class $[X_v]$ in the first Lie algebroid cohomology of $
\Gamma_\beta$, called the \emph{modular class} of $\beta$:
\[
[X_v]\in H^1(M,\Gamma_\beta).
\]

\subsection{Tensor product of Dirac structures}\label{tendir}
We alluded in section~\ref{gerb} to a Baer sum operation on Courant
algebroids; in this section we elaborate on the idea, and introduce
an associated tensor product operation on Dirac structures, which
will be used in Section~\ref{poismod}.
 This operation was noticed independently by the authors of \cite{ABM}, who use it
 to describe some remarkable properties of Dirac structures on Lie
 groups.

Like gerbes, exact Courant algebroids may be pulled back to
submanifolds $\iota:S\hookrightarrow M$.  Following~\cite{BCG}, we
provide a proof in the Appendix.  It is shown there that if $E$ is
an exact Courant algebroid on $M$, then
\[
\iota^*E:=K^\perp/K,
\]
where $K=\Ann(TS)$ and $K^\perp$ is the orthogonal complement in
$E$, inherits an exact Courant algebroid structure over $S$, with
\v{S}evera class given simply by the pullback along the inclusion.
Furthermore, any Dirac structure $L\subset E$ may be pulled back to
$S$ via
\begin{equation}\label{diracpull}
\iota^*L_S:=\frac{L\cap K^\perp+K}{K}\subset \iota^*E,
\end{equation}
which is an integrable Dirac structure whenever it is smooth as a
subbundle of $\iota^*E$, e.g. if $L\cap K^\perp$ has constant rank
on $S$ (see the Appendix, Proposition~\ref{pulldirac}, for a proof).

Let $E,F$ be exact Courant algebroids over the same manifold $M$.
Then $E\times F$ is naturally an exact Courant algebroid over
$M\times M$, and may be pulled back by the diagonal embedding
$d:M\hookrightarrow M\times M$.  The result coincides with the Baer
sum of $E$ and $F$ as defined in~(\cite{Sevlet},\cite{Bressler}),
and we denote it as follows.
\begin{defn}
Let $E_1,E_2$ be exact Courant algebroids over the same manifold
$M$, and let $d:M\lra M\times M$ be the diagonal embedding.  Then we
define the \emph{Baer sum} or \emph{tensor product} of $E_1$ with
$E_2$ to be the exact Courant algebroid (over $M$)
\[ E_1\boxtimes E_2 = d^*(E_1\times E_2),
\]
which can be written simply as
\[
E_1\boxtimes E_2 = \{(e_1,e_2)\in E_1\times E_2\ :\
\pi_1(e_1)=\pi_2(e_2)\}/\{(-\pi_1^*\xi,\pi_2^*\xi)\ :\ \xi\in T^*\},
\]
and has \v{S}evera class equal to the sum $[H_1]+[H_2]$.
\end{defn}
The standard Courant algebroid $(T\oplus T^*,[\cdot,\cdot]_0)$ acts
as an identity element for this operation, and every exact Courant
algebroid $E$ has a natural inverse, denoted by $E^\top$, defined as
the same Courant algebroid with $\IP{\cdot,\cdot}$ replaced with its
negative $-\IP{\cdot,\cdot}$:
\begin{equation}\label{invco}
E^\top = (E,[\cdot,\cdot],-\IP{\cdot,\cdot},\pi).
\end{equation}
Note that this sign reversal changes the sign of $\pi^*:T^*\lra E$
and hence of the curvature $H$ of any splitting~\eqref{splitcurv},
and finally therefore of the \v{S}evera class.
\begin{prop}
Let $E^\top$ be as above.  Then we have a canonical isomorphism
\[
E^\top\boxtimes E = (T\oplus T^*,[\cdot,\cdot]_0).
\]
\end{prop}
\begin{proof}
$E^\top\boxtimes E$ has a well-defined, bracket-preserving splitting
$s:T\lra E^\top\boxtimes E$ given by $X\mapsto [(e_X,e_X)]$ for any
$e_X\in E$ such that $\pi(e_X)=X$.
\end{proof}

We may now use the Dirac pullback~\eqref{diracpull} to define the
tensor product of Dirac structures; an equivalent definition appears
in~\cite{ABM}.
\begin{defn}
Let $L_1\subset E_1$, $L_2\subset E_2$ be Dirac structures and
$E_1,E_2$ as above. We define the \emph{tensor product}
\[
L_1\boxtimes L_2 = d^*(L_1\times L_2)\subset E_1\boxtimes E_2,
\]
where $d^*$ denotes the Dirac pullback~\eqref{diracpull} by the
diagonal embedding.  Explicitly, we have
\begin{equation}\label{dirtens}
L_1\boxtimes L_2=(\{(x_1,x_2)\in L_1\times L_2\ :\
\pi_1(x_1)=\pi_2(x_2)\} + K)/K,
\end{equation}
where $K = \{(-\pi_1^*\xi,\pi_2^*\xi)\ :\ \xi\in T^*\}$. This forms
a Dirac structure whenever it is smooth as a bundle.
\end{defn}

\begin{example}
The canonical Dirac structure $T^*\subset E$ acts as a zero element:
for any other Dirac structure $L\subset F$,
\[
T^*\boxtimes L = T^*\subset E\boxtimes F.
\]
\end{example}
\begin{example}
The Dirac structure $\Delta+\Ann(\Delta)\subset T\oplus T^*$
associated to an integrable distribution $\Delta\subset T$ is
idempotent:
\[
(\Delta + \Ann(\Delta))\boxtimes(\Delta +\Ann(\Delta)) = \Delta
+\Ann(\Delta).
\]
\end{example}
\begin{example}
The tensor product of Dirac structures is compatible with B-field
transformations:
\[
e^{B_1}L_1\boxtimes e^{B_2}L_2 = e^{B_1+B_2}(L_1\boxtimes L_2).
\]
\end{example}
Combining this with the previous example, taking $\Delta=T$, we see
that Dirac structures transverse to $T^*$ remain so after tensor
product.  Finally we provide an example where smoothness is not
guaranteed.
\begin{example}
Let $L\subset E$ be any Dirac structure, with $L^\top\subset E^\top$
defined by the inclusion $L\subset E$.  Then
\[
L^\top\boxtimes L = \Delta +\Ann(\Delta) \subset T\oplus T^*,
\]
where $\Delta = \pi(L)$.  Hence $L^\top \boxtimes L$ is a Dirac
structure when $\Delta+\Ann(\Delta)$ is a smooth subbundle, i.e.
when $\Delta$ has constant rank.
\end{example}

Assuming we choose splittings for $E_1,E_2$, the tensor product of
Dirac structures $L_1\subset E_1,\ L_2\subset E_2$ annihilates the
\emph{wedge product} $K_1\wedge K_2$ of the pure spinor lines
representing $L_1$ and $L_2$.  For reasons of skew-symmetry,
$K_1\wedge K_2$ is nonzero only when $L_1\cap L_2\cap T^*=\{0\}$.
This result also appears in~\cite{ABM}:
\begin{prop}\label{formtens}
Let $L_1,\ L_2$ be Dirac structures in $T\oplus T^*$, and let
$\varphi_1\in K_1$, $\varphi_2\in K_2$ be (local) generators for
their corresponding pure spinor lines in $\wedge^\bullet T^*$.  Then
\[
L_1\boxtimes L_2 \cdot (\varphi_1\wedge \varphi_2) = 0,
\]
and therefore $\varphi_1\wedge\varphi_2$ is a pure spinor for
$L_1\boxtimes L_2$ as long as $L_1\cap L_2\cap T^*=\{0\}$.
\end{prop}
\begin{proof}
From expression~\eqref{dirtens}, we obtain the following simple
expression:
\begin{equation}\label{tensform}
L_1\boxtimes L_2 = \{X+\xi+\eta\ :\ X+\xi\in L_1\ \text{and}\
X+\eta\in L_2\}.
\end{equation}
Then for $X+\xi+\eta\in L_1\boxtimes L_2$, we have
\[
(X+\xi+\eta)\cdot(\varphi_1\wedge\varphi_2) = (i_X\varphi_1 +
\xi\wedge \varphi_1)\wedge \varphi_2 +
(-1)^{|\varphi_1|}\varphi_1\wedge(i_X\varphi_2+\eta\wedge\varphi_2)
= 0.
\]
\end{proof}
There is an anti-orthogonal map $T\oplus T^*\lra T\oplus T^*$ given
by
\[
X+\xi\mapsto (X+\xi)^\top = X - \xi
\]
which satisfies
\[
[(X+\xi)^\top,(Y+\eta)^\top]^\top_H = [X+\xi,Y+\eta]_{-H},
\]
so that it takes the Courant algebroid to its inverse~\eqref{invco}.
This operation intertwines with the Clifford reversal, in the sense
that
\[
((X+\xi)\cdot\varphi)^\top =
(-1)^{|\varphi|+1}(X+\xi)^\top\cdot\varphi^\top,
\]
for any $\varphi\in\wedge^\bullet T^*$, where $|\varphi|$ denotes
the degree.  As a result, we see that reversal operation on forms
corresponds to the reversal $L\mapsto L^\top$ of Dirac structures in
$T\oplus T^*$.  Since the Mukai pairing of pure spinors
$\varphi,\psi$, is given by the top degree component of
$\varphi^\top\wedge \psi$, we conclude from
Propositions~\ref{chevalley} and~\ref{bothways} that for transverse
Dirac structures $L_1,L_2\subset E$ the tensor product
$L_1^\top\boxtimes L_2\subset T\oplus T^*$ has zero intersection
with $T$ and hence is the graph of a Poisson bivector $\beta$.  This
result was first observed in its general form in~\cite{ABM}, and is
consistent with the appearance of a Poisson structure associated to
any Lie bialgebroid in~\cite{MackXu}.
\begin{prop}[Alekseev-Bursztyn-Meinrenken \cite{ABM}]\label{transdir}
Let $E$ be any exact Courant algebroid and $L_1,L_2\subset E$ be
transverse Dirac structures. Then
\[
L_1^\top\boxtimes L_2 = \Gamma_\beta\subset T\oplus T^*,
\]
where $\beta\in C^\infty(\wedge^2 T)$ is a Poisson structure.
\end{prop}

\section{Generalized complex structures}\label{gcg}

Just as a complex structure may be defined as an endomorphism
$J:T\lra T$ satisfying $J^2=-1$ and which is integrable with respect
to the Lie bracket, we have the following definition, due to
Hitchin~\cite{Hitchin}:
\begin{defn}
A \emph{generalized complex structure} on an exact Courant algebroid
$E\cong T\oplus T^*$ is an endomorphism $\JJ:E\lra E$ satisfying
$\JJ^2=-1$ and which is integrable with respect to the Courant
bracket, i.e. its $+i$ eigenbundle $L\subset E\otimes\CC$ is
involutive.
\end{defn}
An immediate consequence of Proposition~\ref{isot} is that the $+i$
eigenbundle of a generalized complex structure must be isotropic,
implying that $\JJ$ must be orthogonal with respect to the natural
pairing on $E$:
\begin{prop}
A generalized complex structure $\JJ$ must be orthogonal, and hence
defines a symplectic structure $\IP{\JJ\cdot,\cdot}$ on $E$.
\end{prop}
\begin{proof}
Let $x,y\in C^\infty(E)$ and decompose $x=a+\bar{a},\ y=b+\bar b$
according to the polarization $E\otimes\CC = L\oplus\overline{L}$.
Since $L$ must be isotropic by Proposition~\ref{isot},
\[
\IP{\JJ x,\JJ y} = \IP{a,\bar b} + \IP{\bar a, b} = \IP{x,y}.
\]
Hence $\JJ$ is orthogonal and $\IP{\JJ\cdot,\cdot}$ is symplectic,
as required.
\end{proof}
This equivalence between complex and symplectic structures on $E$
compatible with the inner product is illustrated most clearly by
examining two extremal cases of generalized complex structures on
$T\oplus T^*$. First, consider the endomorphism of $T\oplus T^*$:
\begin{equation}\label{jcx}
\JJ_J=\left(\begin{matrix}-J&0\\0&J^*\end{matrix}\right),
\end{equation}
where $J$ is a usual complex structure on $V$. Then we see that
$\JJ_J^2=-1$ and $\JJ_J^*=-\JJ_J$.  Its $+i$ eigenbundle $L_J =
T_{0,1}\oplus T^*_{1,0}$ is, by Example~\ref{cxdirac}, integrable if
and only if $J$ is integrable and $H^{(3,0)}=0$.

At the other extreme, consider the endomorphism
\begin{equation}\label{jsymp}
\JJ_\omega=\left(\begin{matrix}0&-\omega^{-1}\\\omega&0\end{matrix}\right),
\end{equation}
where $\omega$ is a usual symplectic structure.  Again, we observe
that $\JJ_\omega^2=-1$ and the $+i$ eigenbundle
\[
L_\omega = \{X - i\omega (X)\ :\ X\in T\otimes\CC\}
\]
is integrable, by Example~\ref{symp}, if and only if $H=0$ and
$d\omega = 0$.

Therefore we see that diagonal and anti-diagonal generalized complex
structures correspond to complex and symplectic structures,
respectively.  We now make some elementary observations concerning
the general case.

\begin{prop}
Generalized complex manifolds must be even-dimensional.
\end{prop}
\begin{proof}
Let $p\in M$ be any point and $E_p$ the fibre of the exact Courant
algebroid at $p$.  Let $x\in E_p$ be null, i.e. $\IP{x,x}=0$.  Then
$\JJ x$ is also null and is orthogonal to $x$.  Therefore $\{x,\JJ
x\}$ span an isotropic subspace $N\subset E_p$. We may iteratively
enlarge the spanning set by adding a pair  $\{x', \JJ x'\}$ for
$x'\in N^\bot$, until $N^\bot=N$ and $\dim M=\dim N$ is even.
\end{proof}

At any point $p\in M$, the orthogonal group $O(E_p)\cong O(2n,2n)$
acts transitively on the space of generalized complex structures at
$p$ by conjugation, with stabilizer $U(n,n)=O(2n,2n)\cap
GL(2n,\CC)$. Therefore the space of generalized complex structures
at $p$ is given by the coset space
\begin{equation}\label{cos}
\frac{O(2n,2n)}{U(n,n)}.
\end{equation}
In this sense, a generalized complex structure on an
even-dimensional manifold is an integrable reduction of the
structure group of $E$ from $O(2n,2n)$ to $U(n,n)$.  Since $U(n,n)$
is homotopic to $U(n)\times U(n)$, the $U(n,n)$ structure may be
further reduced to $U(n)\times U(n)$, which corresponds
geometrically to the choice of a positive definite subbundle
$C_+\subset E$ which is complex with respect to $\JJ$. The
orthogonal complement $C_-=C_+^\bot$ is negative-definite and also
complex, and so we obtain the orthogonal decomposition
\begin{equation}\label{posneg}
E=C_+\oplus C_-.
\end{equation}
Note that since $C_\pm$ are definite and $T^*\subset E$ is
isotropic, the projection $\pi:C_\pm\rightarrow T$ is an
isomorphism.  Hence we can transport the complex structures on
$C_\pm$ to $T$, obtaining two almost complex structures $J_+, J_-$
on $T$.  Thus we see that a generalized complex manifold must admit
an almost complex structure.  Furthermore it has two canonically
associated sets of Chern classes $c^\pm_i=c_i(T,J_\pm)\in
H^{2i}(M,\ZZ)$. Summarizing, and using~\eqref{posneg}, we obtain the
following.
\begin{prop}A generalized complex manifold must admit almost complex
structures, and has two sets of canonical classes $c^\pm_i\in
H^{2i}(M,\ZZ)$ such that the total Chern class
\[
c(E,\JJ) = c^+\cup c^-,
\]
where $c^\pm = \sum_i c^\pm_i$.
\end{prop}

\subsection{Type and the canonical line bundle}\label{typecan}

Any exact Courant algebroid has a canonical Dirac structure
$T^*\subset E$, and a generalized complex structure $\JJ$ may be
characterized by its action on this Dirac structure, as we now
describe.

If $\JJ T^*=T^*$, then $\JJ$ determines a usual complex structure on
the manifold, and a splitting may be chosen for $E$ so that $\JJ$ is
of the form~\eqref{jcx}. On the other hand, if $\JJ T^*\cap
T^*=\{0\}$, then we have the canonical splitting $E=\JJ T^*\oplus
T^*$, and $\JJ$ takes the form~\eqref{jsymp}, i.e. a symplectic
structure.

In general, the subbundle $\JJ T^*\subset E$ projects to a
distribution
\begin{equation}\label{dlta}
\Delta = \pi(\JJ T^*)\subset T
\end{equation}
which may vary in dimension along the manifold.  Defining
\[
E_\Delta = \frac{T^*+\JJ T^*}{\Ann (\Delta)},
\]
we see that $E_\Delta$ is an extension of the form
\[
\xymatrix{0\ar[r]&\Delta^*\ar[r]&E_\Delta\ar[r]^\pi&\Delta\ar[r]&0},
\]
and since $\Ann(\Delta) = T^*\cap\JJ T^*$ is complex, we see that
$\JJ$ induces a complex structure on $E_\Delta$ such that
$\JJ\Delta^*\cap\Delta^* = \{0\}$.  Therefore, at each point,
$\Delta$ inherits a generalized complex structure of symplectic
type. Furthermore,
\[
\frac{E}{T^*+\JJ T^*} = T/\Delta,
\]
showing that, at each point, $T/\Delta$ inherits a complex
structure.  Ignoring integrability, which we address in the next
section, we may therefore conclude that a generalized complex
manifold carries a canonical symplectic distribution (of variable
dimension) with transverse complex structure.

The invariant of $\JJ$ measuring the number of transverse complex
directions at each point is called the \emph{type} of the
generalized complex structure, and may range from $0$, in the case
of a symplectic structure, to $n=\tfrac{1}{2}\dim_\RR M$ for a
complex structure.
\begin{defn}
The \emph{type} of the generalized complex structure $\JJ$ is the
upper semi-continuous function
\[
\mathrm{type}(\JJ)= \tfrac{1}{2}\dim_\RR T^*\cap\JJ T^*,
\]
with possible values $\{0,1,\ldots,n\}$, where
$n=\tfrac{1}{2}\dim_\RR M$.
\end{defn}
The terminology is chosen to coincide with the notion of type for
Dirac structures (see Definition~\ref{type}), since it is indeed the
type of the Dirac structure $L\subset E\otimes\CC$ defining $\JJ$:
\begin{prop}
The type of $\JJ$ coincides with the type of its $+i$ eigenbundle
$L\subset E\otimes\CC$, and hence is of fixed parity throughout the
manifold.
\end{prop}
\begin{proof}
At any point, the subspace $T^*\cap\JJ T^*$ is complex, and hence
\[
(T^*\cap \JJ T^*)\otimes\CC = A\oplus \overline{A},
\]
where $A=L\cap(T^*\otimes\CC)$. Since $\mathrm{type}(L)=\dim_\CC
L\cap(T^*\otimes\CC)$, we see that
$\mathrm{type}(\JJ)=\mathrm{type}(L)$, as required.  By
Proposition~\ref{partype}, the parity of $\type(\JJ)$ must be fixed
throughout the manifold.
\end{proof}
As a result, we see that in real dimension $2$, connected
generalized complex manifolds must be of constant type $0$ or $1$,
i.e. of symplectic or complex type, whereas in dimension 4, they may
be of types $0$, $1$, or $2$, with possible jumping from $0$
(symplectic) to $2$ (complex) along closed subsets of the manifold;
we shall encounter such examples in sections~\ref{examples}
and~\ref{deformex}.

It also follows from our work on Dirac structures that a generalized
complex structure is completely characterized by the pure spinor
line $K\subset S\otimes\CC$ corresponding to the maximal isotropic
subbundle $L$.  When a splitting for $E$ is chosen, we obtain an
identification $S=\wedge^\bullet T^*\otimes(\det T)^{1/2}$, and
hence $K$ may be viewed as a line subbundle of the complex
differential forms. For a symplectic structure,
$L_\omega=e^{-i\omega}(T)$, and so
\[
K_\omega = e^{-i\omega}\cdot \wedge^0 T^*  = \CC \cdot e^{i\omega},
\]
whereas for a complex structure $L_J = T_{0,1}+T^*_{1,0}$, so that
\[
K_J = \wedge^n T^*_{1,0},
\]
leading to the following definition.
\begin{defn}
The \emph{canonical line bundle} of a generalized complex structure
on
 $T\oplus T^*$ is the complex pure spinor line subbundle $K\subset
\wedge^\bullet T^*\otimes\CC$ annihilated by the $+i$ eigenbundle
$L$ of $\JJ$.
\end{defn}
Proposition~\ref{evenodd} states that a generator $\varphi\in K_x$
for the canonical line bundle at the point $x\in M$ must have the
form
\begin{equation}\label{genr}
\varphi = e^{B+i\omega}\Omega,
\end{equation}
where $\Omega = \theta_1\wedge\cdots\wedge\theta_k$ for
$(\theta_1,\ldots,\theta_k)$ a basis for $L\cap(T^*\otimes\CC)$, and
$B,\omega$ are the real and imaginary components of a complex
2-form. As a result we can read off the type of $\JJ$ at $p$
directly as the least nonzero degree ($k$) of the differential form
$\varphi$. The generalized complex structure defines a polarization
\begin{equation}
E\otimes\CC = L\oplus\overline{L},
\end{equation}
and therefore by Proposition~\ref{chevalley},
\begin{equation}\label{nond}
\IP{\varphi,\overline\varphi}\neq 0.
\end{equation}
Using~\eqref{genr}, we obtain
\begin{align*}
0\neq\IPS{(e^{B+i\omega}\Omega,e^{B-i\omega}\overline\Omega)}&=\IPS{(e^{2i\omega}\Omega,\overline\Omega)}\\
&=\tfrac{(-1)^{2n-k}(2i)^{n-k}}{(n-k)!}\omega^{n-k}\wedge
\Omega\wedge\overline\Omega,
\end{align*}
which expresses the fact that $\omega$ pulls back to the symplectic
form on $\Delta=\ker\Omega\wedge\overline\Omega$ described earlier,
and $\Omega$ defines the complex structure transverse to $\Delta$.
We also see that $\IP{\varphi,\overline\varphi}\in\det T^*$ defines
an orientation independent of the choice of $\varphi$, giving a
global orientation on the manifold. This orientation, together with
the parity of the type, defines a pair of invariants which
distinguish the four connected components of the coset
space~\eqref{cos}.

In the following result, we show that at any point, a splitting for
$E$ may be chosen so that the generalized complex structure is a
product of a complex and a symplectic structure of lesser dimension.
\begin{theorem}
At any point, a generalized complex structure of type $k$ is
equivalent, by a choice of splitting for $E$, to the direct sum of a
complex structure of complex dimension $k$ and a symplectic
structure of real dimension $2n-2k$.
\end{theorem}
\begin{proof}
Fixing a splitting for $E$ at $x\in M$, the generalized complex
structure is defined, as in~\eqref{genr}, by the pure spinor
\begin{equation*}
\varphi=e^{B+i\omega}\Omega,
\end{equation*}
where $\omega^{n-k}\wedge\Omega\wedge\overline{\Omega}\neq 0$.
Choose a subspace $N\subset T_x$ transverse to
$\Delta=\ker\Omega\wedge\overline{\Omega}$.  Then $\Delta$ carries a
symplectic structure $\omega_0=\omega|_{\Delta}$ and $N$ inherits a
complex structure determined by $\Omega|_N$. The 2-forms then
decompose as
\begin{equation*}
\wedge^2
T^*_x=\bigoplus_{p+q+r=2}\wedge^p\Delta^*\otimes\wedge^qN^*_{1,0}\otimes\wedge^rN^*_{0,1},
\end{equation*}
so that forms have tri-degree $(p,q,r)$. While $\Omega$ is purely of
type $(0,k,0)$, the complex 2-form $A=B+i\omega$ decomposes into six
components:
\begin{equation*}\begin{array}{ccc}
 A^{200}  &  &  \\
  A^{110} & A^{101} &  \\
  A^{020} & A^{011} & A^{002}
\end{array}\end{equation*}
Only the components $A^{200},A^{101},A^{002}$ act nontrivially on
$\Omega$ in the expression $e^A\Omega$.  Hence we are free to modify
the other three components at will.  Note that
$\omega_0=-\tfrac{i}{2}(A^{200}-\overline{A^{200}})$. Now define the
real 2-form
\begin{equation*}
\tilde
B=\tfrac{1}{2}(A^{200}+\overline{A^{200}})+A^{101}+\overline{A^{101}}+A^{002}+\overline{A^{002}},
\end{equation*}
and observe that $e^{\tilde B +
i\omega_0}\Omega=e^{B+i\omega}\Omega$, demonstrating that
$\varphi=e^{\tilde B+i\omega_0}\Omega$, i.e. $\varphi$ is a B-field
transform of $e^{i\omega_0}\Omega$, which is a direct sum of a
symplectic structure on $\Delta$ and complex structure on $N$, as
required.
\end{proof}

The canonical line bundle $K$ introduced in this section, along with
its complex conjugate $\overline K$, are the extremal line bundles
of a $\ZZ$-grading on spinors induced by the generalized complex
structure.  As described in~\eqref{altgrad}, a polarization induces
a $\ZZ$-grading on spinors;  therefore a generalized complex
structure on $T\oplus T^*$, since it determines a polarization
$(T\oplus T^*)\otimes\CC = L\oplus \overline{L}$, induces an
alternative $\ZZ$-grading on differential forms
\[
\wedge^\bullet T^*\otimes\CC = U^{-n}\oplus\cdots\oplus U^n,
\]
where $U^{n}=K$ is the canonical line bundle and $U^{n-k} = \wedge^k
\overline{L}\cdot U^{n}$.  Since $\overline L$ annihilates $U^{-n}$,
we see that $U^{-n} = \overline{U^{n}}$ is the canonical line of
$-\JJ$.  We therefore have the following convenient description of
this $\ZZ$-grading.

\begin{prop}A generalized complex structure $\JJ$ on $E=T\oplus T^*$ gives rise to a $\ZZ$-grading
\[
\wedge^\bullet T^*\otimes\CC = U^{-n}\oplus\cdots\oplus U^n,
\]
where $U^k$ is the $ik$-eigenbundle of $\JJ$ acting in the spin
representation, and $U^n=K$, the canonical line bundle.
\end{prop}
In the case of a usual complex structure $\JJ_J$, then the graded
components correspond to the well-known $(p,q)$-decomposition of
forms as follows:
\begin{equation}\label{pqform}
U^k_J = \bigoplus_{p-q=k} \Omega^{p,q}(M,\CC),
\end{equation}
since $\JJ_J$ acts via the spin representation as $J^*$, which has
eigenvalue $i(p-q)$ on $\Omega^{p,q}$.

The fact that $U^{-n}=\det\overline{L}\cdot U^{n}$, combined with
our previous remark~\eqref{nond}, implies that
\[
U^{n}\otimes\det L^*\otimes U^n \cong \det T^*\otimes\CC.
\]
Since the complex bundle $L$ is isomorphic to $(E,\JJ)$, we obtain
the following.
\begin{corollary}
The canonical line bundle of a generalized complex manifold has
first Chern class satisfying
\[
2c_1(K) ={c_1^+ + c_1^-}.
\]
\end{corollary}

\subsection{Courant integrability}

The notion of type and the $\ZZ$-grading on spinors introduced in
the last section do not depend on the Courant integrability of the
generalized complex structure; they may be associated to any
\emph{generalized almost complex} structure:
\begin{defn}
A \emph{generalized almost complex} structure is a complex structure
$\JJ$ on an exact Courant algebroid which is orthogonal with respect
to the natural inner product.
\end{defn}
Naturally, a generalized almost complex structure $\JJ$ is said to
be integrable when its $+i$ eigenbundle $L\subset E\otimes\CC$ is
involutive for the Courant bracket, i.e. $L$ is a Dirac structure.
\begin{prop}
A generalized complex structure is equivalent to a complex Dirac
structure $L\subset E\otimes\CC$ such that
$L\cap\overline{L}=\{0\}$.
\end{prop}
As a result, $(L,[,],\pi)$, where $\pi:L\lra T\otimes\CC$ is the
projection, defines the structure of a Lie algebroid, and therefore
we obtain a differential complex
\begin{equation}\label{ellcx}
\xymatrix{C^\infty(\wedge^k L^*)\ar[r]^{d_L}&
C^\infty(\wedge^{k+1}L^*)},
\end{equation}
where $d_L$ is the Lie algebroid de Rham differential, which
satisfies $d_L^2=0$ due to the Jacobi identity for the Courant
bracket restricted to $L$. The operator $d_L$ has principal symbol
$s(d_L):T^*\otimes\wedge^k L^*\rightarrow \wedge^{k+1}L^*$ given by
$\pi^*:T^*\rightarrow L^*$ composed with wedge product, i.e.
\begin{equation*}
s_\xi(d_L)=\pi^*(\xi)\wedge\cdot\ \ ,
\end{equation*}
where $\xi\in T^*$.  We now observe that the complex~\eqref{ellcx}
is \emph{elliptic} for a generalized complex structure.
\begin{prop}\label{ellipcx}
The Lie algebroid complex of a generalized complex structure is
elliptic.
\end{prop}
\begin{proof}
Given a real, nonzero covector $\xi\in T^*$, write $\xi = \alpha
+\overline\alpha$ for $\alpha\in L$.  For $v\in L$, we have
$\pi^*\xi(v) = \xi(\pi(v)) = \IP{\xi,v} = \IP{\overline\alpha,v}$.
Using the inner product to identify $L^*=\overline{L}$, we therefore
have $\pi^*\xi = \overline{\alpha}$, which is clearly nonzero if and
only if $\xi$ is.   As a result, the symbol sequence is exact for
any nonzero real covector, as required.
\end{proof}
This provides us with our first invariants associated to a
generalized complex structure:
\begin{corollary}
The cohomology of the complex~\eqref{ellcx}, called the Lie
algebroid cohomology $H^\bullet(M,L)$, is a finite dimensional
graded ring associated to any compact generalized complex manifold.
\end{corollary}
In the case of a complex structure, $L=T_{0,1}\oplus T^*_{1,0}$,
while $d_L = \delbar$, and so the Lie algebroid complex is a sum of
usual Dolbeault complexes, yielding
\[
H^k(M,L_J) = \bigoplus_{p+q=k} H^p(M,\wedge^qT_{1,0}).
\]
In the case of a symplectic structure, the Lie algebroid $L$ is the
graph of $i\omega$, and hence is isomorphic to $T\otimes\CC$ as a
Lie algebroid.  Hence its Lie algebroid cohomology is simply the
complex de Rham cohomology.
\[
H^k(M,L_\omega)= H^k(M,\CC).
\]

We now describe a second invariant, obtained from the $\ZZ$-grading
on differential forms induced by $\JJ$.  As we saw in the previous
section, a generalized complex structure on $T\oplus T^*$ determines
an alternative grading for the differential forms
\begin{equation*}
\wedge^\bullet T^*\otimes\CC=U^{-n}\oplus \cdots\oplus U^{n}.
\end{equation*}
This $\ZZ$-grading may be viewed as the intersection of two complex
conjugate filtrations
\begin{equation*}
F_i = \oplus_{k=0}^{i} U^{n-k},\ \ \ \ \overline{F_i} =
\oplus_{k=0}^{i} U^{-n+k}.
\end{equation*}
More precisely, we have
\begin{equation}\label{intfilt}
U^k = F_{n-k}\cap \overline F_{n+k}.
\end{equation}
By Theorem~\ref{dfilt}, the integrability of $\JJ$ with respect to
$\Cour{\cdot,\cdot}_H$ is equivalent to the fact that $d_H$ takes
$C^\infty(F_i)$ into $C^\infty(F_{i+1})$.  Using~\eqref{intfilt},
this happens if and only if $d_H$ takes $C^\infty(U^k)$ into
$C^\infty(U^{k-1}\oplus U^k\oplus U^{k+1})$, but since $d_H$ is odd,
 we see that $\JJ$ is integrable if and only if $d_H(C^\infty(U^k))\subset C^\infty(U^{k-1}\oplus
U^{k+1})$.  Projecting to these two components, we obtain
$d_H=\del+\delbar$, where
\begin{equation}\label{gendol}
\xymatrix{C^\infty(U^k)\ar@<0.5ex>[r]^\del&C^\infty(U^{k+1})\ar@<0.5ex>[l]^\delbar}.
\end{equation}
In greater generality, we may use our calculation in~\eqref{tcour}
to obtain the following.
\begin{theorem}\label{deldelbar}
Let $\JJ$ be a generalized almost complex structure on $T\oplus
T^*$, and define
\begin{align*}
\del&=\pi_{k+1}\circ d_H:C^\infty(U^k)\lra C^\infty
(U^{k+1})\\
\delbar&=\pi_{k-1}\circ d_H:C^\infty(U^k)\lra C^\infty (U^{k-1}),
\end{align*}
where $\pi_k$ is the projection onto $U^k$.  Then
\begin{equation}\label{noint}
d_H = \del + \delbar + T_L + \overline{T_L},
\end{equation}
where $T_L\in\wedge^3 L^*=\wedge^3\overline{L}$ is defined by
\[
T_L(e_1,e_2,e_3)=\IP{[e_1,e_2],e_3},
\]
and acts via the Clifford action in~\eqref{noint}.  $\JJ$ is
integrable, therefore, if and only if $d_H=\del+\delbar$, or
equivalently, if and only if
\begin{equation}
d_H(C^\infty(U^n))\subset C^\infty(U^{n-1}).
\end{equation}
\end{theorem}

In the integrable case, since $d_H=\del+\delbar$ and $d_H^2=0$, we
conclude that $\del^2 = \delbar^2 = 0$ and $\del\delbar
=-\delbar\del$; hence in each direction,~\eqref{gendol} defines a
differential complex.

\begin{remark}
Given the above, a generalized complex structure gives rise to a
real differential operator $d^{\JJ}=i(\delbar-\del)$, which can also
be written $d^{\JJ}=[d,\JJ]$, and which satisfies $(d^{\JJ})^2=0$.
It is interesting to note that while in the complex case $d^{\JJ}$
is just the usual $d^c$-operator $d^c=i(\delbar-\del)$, in the
symplectic case $d^\JJ$ is equal to the symplectic adjoint of $d$
defined by Koszul~\cite{Koszul} and studied by
Brylinski~\cite{Brylinski} in the context of symplectic harmonic
forms.
\end{remark}
Using the identification $U^{n-k} = \wedge^kL^*\otimes K$ as
in~\eqref{altgrad}, the operator $\delbar$ can be viewed as a Lie
algebroid connection
\begin{equation*}
\delbar: C^\infty(\wedge^k L^*\otimes K)\rightarrow
C^\infty(\wedge^{k+1}L^*\otimes K),
\end{equation*}
extended from $d_H:C^\infty(K)\rightarrow C^\infty(L^*\otimes K)$
via the rule
\begin{equation}\label{leib}
\delbar(\mu\otimes s)=d_L\mu\otimes s + (-1)^{|\mu|}\mu\wedge ds,
\end{equation}
for $\mu\in C^\infty(\wedge^k L^*)$ and $s\in C^\infty(K)$, and
satisfying $\delbar^2=0$. Therefore $K$ is a module for the Lie
algebroid $L$, and we may call it a \emph{generalized holomorphic
bundle}. From the ellipticity of the Lie algebroid complex for $L$
and the fact that $K$ is a module over $L$, we immediately obtain
the following.
\begin{prop}
The cohomology of the complex $(U^\bullet,\delbar)$, called the
\emph{generalized Dolbeault cohomology} $H^\bullet_{\delbar}(M)$, is
a finite dimensional graded module over $H^\bullet(M,L)$ associated
to any compact generalized complex manifold.
\end{prop}
In the case of a complex structure, Equation~\eqref{pqform} shows
that the generalized Dolbeault cohomology coincides with the usual
Dolbeault cohomology, with grading
\[
H^k_{\delbar}(M) = \bigoplus_{p-q=k} H^{p,q}_{\delbar}(M).
\]
A special case occurs when the canonical line bundle is
\emph{holomorphically trivial}, in the sense that $(K,\delbar)$ is
isomorphic to the trivial bundle $M\times\CC$ together with the
canonical Lie algebroid connection $d_L$.  Then the Lie algebroid
complex and the generalized Dolbeault complex $(U^\bullet,\delbar)$
are isomorphic and hence $H^\bullet_{\delbar}(M)\cong
H^\bullet(M,L)$. This holomorphic triviality of $K$ is equivalent to
the existence of a nowhere-vanishing section $\rho\in C^\infty(K)$
satisfying $d_H\rho=0$. In~\cite{Hitchin}, Hitchin calls these
generalized Calabi-Yau structures:
\begin{defn}
A \emph{generalized Calabi-Yau structure} is a generalized complex
structure with holomorphically trivial canonical bundle, i.e.
admitting a nowhere-vanishing $d_H$-closed section $\rho\in
C^\infty(K)$.
\end{defn}
An example of a generalized Calabi-Yau structure is of course the
complex structure of a Calabi-Yau manifold, which admits a
holomorphic volume form $\Omega$ trivializing the canonical line
bundle.  On the other hand, a symplectic structure has canonical
line bundle generated by the closed form $e^{i\omega}$, so it too is
generalized Calabi-Yau.

Assuming that the canonical bundle is trivial as a smooth line
bundle, i.e. $c_1(K)=0$, we may always choose a non-vanishing
section $\rho\in C^\infty(K)$; by Theorem~\ref{deldelbar},
integrability implies that
\[
d_H\rho = \chi_\rho\cdot\rho,
\]
for a uniquely determined $\chi_\rho\in C^\infty(\overline
L)=C^\infty(L^*)$.  Applying~\eqref{leib}, we obtain
\[
0=d_H^2\rho = (d_L\chi_\rho)\cdot\rho -
\chi_\rho\cdot(\chi_\rho\cdot\rho),
\]
implying that $d_L\chi_\rho=0$.  Just as for the modular class of a
Poisson structure~\eqref{mvf}, $\chi_\rho$ defines a class in the
Lie algebroid cohomology
\begin{equation}\label{obscy}
[\chi_\rho]\in H^1(M,L)
\end{equation}
which is the obstruction to the existence of generalized Calabi-Yau
structure.

More generally, we may use standard \v{C}ech arguments to show that
any generalized holomorphic line bundle $V$ is classified up to
isomorphism by an element $[V]\in\HH^1(\LL_{\text{log}})$ in the
first hypercohomology of the complex $\LL_{\text{log}}$, given by
\[
\xymatrix{C^\infty(\CC^*)\ar[r]^{d_L\log}&C^\infty(L^*)\ar[r]^{d_L}&C^\infty(\wedge^2
L^*)\ar[r]^{d_L}&\cdots}.
\]
\begin{defn}
The Picard group of isomorphism classes of rank 1 generalized
holomorphic bundles, i.e. modules over $L$, is
$Pic(\JJ)=\HH^1(\LL_{\text{log}})$.
\end{defn}
Of course this implies that $\JJ$ is generalized Calabi-Yau if and
only if $[K]=0$ as a class in $\HH^1(\LL_{\text{log}})$. The usual
exponential map induces a long exact sequence of hypercohomology
groups
\[
\xymatrix{ \cdots\ar[r]& H^1(M,L)\ar[r] &
\HH^1(\LL_{\log})\ar[r]^{c_1} &H^{2}(\ZZ)\ar[r] &\cdots },
\]
and so we recover the observation~\eqref{obscy} that when $c_1(K)=0$
the Calabi-Yau obstruction lies in $H^1(M,L)$.

\begin{example}\label{gcholb}
Suppose that the complex bundle $V$ is generalized holomorphic for a
complex structure $\JJ_J$.  Then the differential $D:C^\infty(V)\lra
C^\infty(L^*\otimes V )$ may be decomposed according to
$L=T_{0,1}\oplus T^*_{1,0}$ to yield
\[
D=\delbar_V+\Phi,
\]
where $\delbar_V:C^\infty(V)\lra C^\infty(T^*_{0,1}\otimes V)$ is a
usual partial connection, $\Phi:V\lra T_{1,0}\otimes V$ is a bundle
map, and $D\circ D=0$ yields the conditions
\begin{itemize}
\item $\delbar_V^2=0$, i.e. $V$ is a usual holomorphic bundle,
\item $\delbar_V(\Phi)=0$, i.e. $\Phi$ is holomorphic,
\item $\Phi\wedge\Phi=0$ in $\wedge^2 T_{1,0}\otimes\End(V)$.
\end{itemize}
In the rank 1 case, therefore, we obtain the result
\[
Pic(\JJ_J) = H^1(\OO^*)\oplus H^0(\mathcal{T}),
\]
showing that the generalized Picard group contains the usual Picard
group of the complex manifold but also includes its infinitesimal
automorphisms.
\end{example}

\subsection{Hamiltonian symmetries}

The Lie algebra $\textbf{sym}(\omega)$ of infinitesimal symmetries
of a symplectic manifold consists of sections $X\in C^\infty(T)$
such that $\LL_X\omega = 0$. The \emph{Hamiltonian} vector fields
$\textbf{ham}(\omega)$ are those infinitesimal symmetries generated
by smooth functions, in the sense $X = \omega^{-1}(df)$, for $f\in
C^\infty(M,\RR)$.  We then have the well-known sequence
\[
\xymatrix{0\ar[r]&\textbf{ham}(\omega)\ar[r]&\textbf{sym}(\omega)\ar[r]^{\omega^{-1}}&H^1(M,\RR)\ar[r]&0}.
\]
We now give an analogous description of the symmetries of a
generalized complex structure and examine the manner in which it
specializes in the cases of symplectic and complex geometry.
\begin{defn}
An infinitesimal symmetry $v\in\text{{\bf sym}}(\JJ)$ of a
generalized complex structure $\JJ$ on the Courant algebroid $E$ is
defined to be a section $v\in C^\infty(E)$ which preserves $\JJ$
under the adjoint action, i.e. $\ad_v\circ\JJ=\JJ\circ\ad_v$, or
equivalently, $[v,C^\infty(L)]\subset C^\infty(L)$.
\end{defn}
In the presence of a generalized complex structure $\JJ$, a real
section $v\in C^\infty(E)$ may be decomposed according to the
splitting $E\otimes\CC = L\oplus \overline{L}$, yielding $v= v^{1,0}
+ v^{0,1}$.  Clearly $[v^{1,0},C^\infty(L)]\subset C^\infty(L)$ by
the integrability of $\JJ$.  However $[v^{0,1},C^\infty(L)]\subset
C^\infty(L)$ if and only if $d_Lv^{0,1}=0$, where we use the
identification $\overline L = L^*$. As a result we identify
$\textbf{sym}(\JJ)= \ker d_L\cap C^\infty(L^*)$, and the
differential complex~\eqref{ellcx} provides the following sequence,
suggesting the definition of generalized Hamiltonian symmetries:
\[
\xymatrix{C^\infty(M,\CC)\ar[r]^{d_L}&\textbf{sym}(\JJ)\ar[r]&H^1(M,L)\ar[r]&0}.
\]
\begin{defn}
An infinitesimal symmetry $v\in\text{\bf sym}(\JJ)$ is
\emph{Hamiltonian}, i.e. $v\in\text{\bf ham}(\JJ)$, when $v = Df$
for $f\in C^\infty(M,\CC)$, where
\[
Df = d_L f + \overline{d_L f} = d(\mathrm{Re} f) - \JJ d(\mathrm{Im}
f).
\]
\end{defn}
As a result we obtain the following exact sequence of complex vector
spaces:
\[
\xymatrix{0\ar[r]&\textbf{ham}(\JJ)\ar[r]&\textbf{sym}(\JJ)\ar[r]&H^1(M,L)\ar[r]&0}.
\]
In the case of a symplectic structure, a section $X+\xi\in
C^\infty(T\oplus T^*)$ preserves $\JJ_{\omega}$ precisely when
$\LL_X\omega=0$ and $d\xi=0$.
 On the other hand, computing $Df$, we obtain
\[
Df =  d(\mathrm{Re}f) + \omega^{-1}d(\mathrm{Im} f),
\]
showing that $X+\xi$ is Hamiltonian precisely when $X$ is
Hamiltonian and $\xi$ is exact.

In the complex case, $X+\xi$ preserves $\JJ_J$ exactly when
$\delbar(X^{1,0}+\xi^{0,1})=0$, i.e. when $X$ is a holomorphic
vector field and $\delbar\xi^{0,1}=0$.   We also have
\[
Df = \delbar f + \del \bar{f},
\]
showing that $X+\xi$ is Hamiltonian exactly when $X=0$ and
$\xi=\delbar f + \del \bar f$ for $f\in C^\infty(M,\CC)$.

Even for a usual complex manifold, therefore, there are nontrivial
Hamiltonian symmetries $\xi=\delbar f + \del\overline{f}$, which
integrate to $B$-field transformations $e^{tB}$, for
$B=\del\delbar(f - \overline{f})$.

\subsection{The Poisson structure and its modular class}\label{poismod}
In this section we describe a natural Poisson
structure on a generalized complex manifold which governs the
behaviour of the symplectic distribution~$\Delta$ introduced in
Section~\ref{typecan}. A formulation of the integrability of $\JJ$
which will be of use is the analog of the vanishing of the Nijenhuis
tensor of an almost complex structure.
\begin{defn}
Let $\JJ$ be a generalized almost complex structure. Then we define
the Nijenhuis tensor $N_\JJ\in C^\infty(\wedge^2 E^*\otimes E)$ as
follows:
\begin{equation}\label{nij}
N_\JJ(e_1,e_2) = [\JJ e_1,\JJ e_2] - \JJ[\JJ e_1, e_2] - \JJ[e_1,\JJ
e_2] - [e_1,e_2].
\end{equation}
\end{defn}
As in the case of an almost complex structure, $\JJ$ is integrable
if and only if $N_\JJ=0$ by the usual argument, which we omit.

The endomorphism $\JJ$ gives rise to an orthogonal $S^1$-action on
the total space of the Courant algebroid $E$; indeed we have, for
$v\in E$ and $t\in\RR$,
\[
e^{it}\cdot v = e^{t\JJ}(v).
\]
We now show that when $\JJ$ is integrable, the $S^1$ family of
almost Dirac structures obtained by applying the above action to
$T^*$ is actually integrable for all $t$.
\begin{prop}\label{circlecx}
Let $\JJ$ be a generalized complex structure. Then the family of
almost Dirac structures
\[
D_t = e^{t\JJ}(T^*)
\]
is integrable for all $t$.
\end{prop}
\begin{proof}
Let $a,b\in\RR$.  Then for $\xi,\eta\in C^\infty(T^*)$, we
use~\eqref{nij} to obtain
\begin{align*}
[(a+b\JJ)\xi,(a+b\JJ)\eta]&=ab([\xi,\JJ\eta] + [\JJ\xi,\eta]) +
b^2[\JJ\xi,\JJ\eta]\\
&=b(a+b\JJ)([\xi,\JJ\eta] + [\JJ\xi,\eta]).
\end{align*}
Since $([\xi,\JJ\eta] + [\JJ\xi,\eta])$ is a 1-form, we see that
$(a+b\JJ)T^*$ is involutive.  Setting $a=\cos t$ and $b=\sin t$ for
each $t$, we obtain the result.
\end{proof}
The path of Dirac structures $D_t$ may be differentiated at $t=0$ as
a path in the Grassmannian of maximal isotropic subbundles of $E$,
yielding a bundle map $P:T^*\lra E/T^*=T$, given by the expression,
for $\xi,\eta\in T^*$,
\[
P(\xi,\eta)=\tfrac{d}{dt}{D_t}(\xi,\eta) =
\tfrac{d}{dt}\IP{e^{t\JJ}\xi,\eta} = \IP{\JJ\xi,\eta}.
\]
As a map $T^*\lra T$, therefore, $P=\pi\circ \JJ$. Therefore we see
immediately that $\im P=\Delta$, defined in Equation~\eqref{dlta}.
We now show that $P$ is a Poisson structure.
\begin{prop}\label{projpois}
The bivector field $P=\pi\circ \JJ|_{T^*}:T^*\lra T$ is Poisson.
\end{prop}
\begin{proof}
Choose a splitting for the Courant algebroid, with curvature $H$.
Then for sufficiently small $t$, the Dirac structures $D_t$ may be
described as graphs of bivector fields $\beta_t:T^*\lra T$.  Since
$D_t$ are integrable, Example~\ref{pois} indicates that $\beta_t$
satisfy
\[
[\beta_t,\beta_t] = \wedge^3\beta_t^*(H).
\]
Letting $t\rightarrow 0$, we see that the cubic term is negligible
and $[P,P]=0$, as required.
\end{proof}
\begin{corollary}\label{foll}
The distribution $\Delta=\pi(\JJ T^*)=\im P$ integrates to a
generalized foliation by smooth symplectic leaves with codimension
$2k$, where $k=\type(\JJ)$.
\end{corollary}
Given an isotropic splitting for the Courant algebroid, $\JJ$ may be
written as a block matrix
\begin{equation}\label{blockm}
\JJ=\left(%
\begin{array}{cc}
  A& P \\
  \sigma & -A^*
\end{array}%
\right),
\end{equation}
so that the Poisson tensor $P$ is apparent.  For a direct
calculation that $P$ is Poisson, as well as more details concerning
its relationship to the tensors $A,\sigma$, see~\cite{Lindstrom};
(see also~\cite{abo},\cite{crainic}).

In the preceding discussion, the Poisson structure $P$ had a natural
 interpretation as an infinitesimal deformation of the Dirac structure $T^*$ rather than
 a genuine Dirac structure.  We now use the tensor product of Dirac
structures described in Section~\ref{tendir} to provide an
alternative, more global, interpretation of $P$ as a Dirac structure
in $T\oplus T^*$. In particular, Proposition~\ref{transdir} suggests
the following result.
\begin{prop}
Let $\JJ$ be a generalized complex structure with $+i$ eigenbundle
$L\subset E\otimes \CC$.  Then
\begin{equation}\label{ptens}
L^\top\boxtimes \overline L = \Gamma_{{iP}/{2}},
\end{equation}
i.e. the tensor product of $L^\top$ with $\overline{L}$ is the graph
of the Poisson structure $iP/2$ in $(T\oplus T^*)\otimes\CC$.
\end{prop}
\begin{proof}
Let $\zeta\in T^*\otimes\CC$, so that $\zeta - i\JJ\zeta\in L$, or
in any splitting, using~\eqref{blockm}, we have $\zeta - iP\zeta +
iA^*\zeta \in L$.  Therefore $(\zeta - iP\zeta + iA^*\zeta)^\top =
(-\zeta - iP\zeta - iA^*\zeta)\in L^\top$ and $\zeta + iP\zeta -
iA^*\zeta \in \overline L$. Combining these using~\eqref{tensform},
we see that
\[
iP\zeta + 2\zeta \in L^\top\boxtimes \overline L,
\]
and hence $\Gamma_{iP/2}\subset L^\top\boxtimes\overline L$. Since
both sides are maximal isotropic subbundles, we must have equality,
as required.
\end{proof}
Besides the fact that this provides an alternative proof of the fact
that $P$ is Poisson, it also relates the Lie algebroids defined by
$L$ and $\overline L$ to that defined by the Poisson structure $P$.
We now observe that this implies a relation between the Calabi-Yau
obstruction class and the modular class.
\begin{prop}
Let $\JJ$ be a generalized complex structure such that $c_1(K)=0$,
and let $\rho\in C^\infty(K)$ be a non-vanishing section with
\[
d_H\rho=v\cdot\rho,\ \ v\in C^\infty(E).
\]
Then $-2\pi(\JJ v)=X$ is the modular vector field associated to the
Poisson structure $P$ and volume form $\IPS{(\rho,\overline\rho)}$.
\end{prop}
\begin{proof} Let $d_H\rho =v^{0,1}\cdot\rho$ for uniquely
defined $v^{0,1}\in C^\infty(\overline L)$, so that $v=v^{1,0} +
v^{0,1}$ for $v^{1,0}=\overline{v^{0,1}}$. By Equation~\eqref{ptens}
and Proposition~\ref{formtens}, we have that
\[
\rho^\top\wedge\overline\rho =
e^{\frac{iP}{2}}\IPS{(\rho,\overline\rho)}=\varphi.
\]
Taking the exterior derivative, and using the definition~\eqref{mvf}
of the modular vector field, we have
\[
\begin{split} d\varphi
=\tilde X\cdot\varphi&= (-1)^{|\rho|}((d_H\rho)^\top\wedge \overline
\rho +
\rho^\top\wedge(d_H\overline\rho) )\\
&=(-1)^{|\rho|}((v^{0,1}\cdot\rho)^\top\wedge\overline\rho
+\rho^\top\wedge(v^{1,0}\cdot\overline\rho))\\
&=-\pi(v^{0,1} - v^{1,0})\cdot(\rho^\top\wedge\overline\rho)\\
&=-i\pi(\JJ v)\cdot \varphi,
\end{split}
\]
showing that $\tilde X=-i\pi(\JJ v)$ is the modular vector field for
$iP/2$.  Rescaling the Poisson structure, we obtain the result.
\end{proof}
\begin{corollary}
The Poisson structure $P$ associated to a generalized Calabi-Yau
manifold is \emph{unimodular} in the sense of
Weinstein~\cite{Weinmod}, i.e. it has vanishing modular class.
\end{corollary}
The map $H^1(M,L)\lra H^1(M,\Gamma_{P})$ of Lie algebroid cohomology
groups implicit in the above result may be understood from the fact
that the projection map $T^*\otimes\CC\lra L$ obtained from the
splitting $E\otimes\CC = L\oplus \overline L$, is actually a Lie
algebroid morphism, when $T^*\otimes\CC$ is endowed with the Poisson
Lie algebroid structure, as we now explain.
\begin{prop}\label{morphismo}
Let $L, P$ be the $+i$-eigenbundle and Poisson structure associated
to a generalized complex structure. The bundle map $a:\Gamma_P\lra
L$ given, for any $\xi\in T^*\otimes\CC$, by
\begin{equation}\label{aye}
a:\xi + P\xi\mapsto i\xi+\JJ \xi,
\end{equation}
is a Lie algebroid homomorphism.
\end{prop}
\begin{proof}
The map $a$ commutes with the projections to the tangent bundle,
since $P=\pi\circ\JJ|_{T^*}$. Given 1-forms $\xi,\eta$, we have
\begin{align*}
[a(\xi+P\xi),a(\eta + P\eta)] &= i([\xi,\JJ\eta] + [\JJ\xi,\eta]) +
[\JJ\xi,\JJ\eta]\\  &=i([\xi,P\eta] + [P\xi,\eta])+\JJ([\xi,P\eta] +
[P\xi,\eta])\\ &= a([\xi+P\xi,\eta+P\eta]),
\end{align*}
as required.
\end{proof}
As a final example of the relationship between a generalized complex
structure and its associated Poisson structure, we use the above Lie
algebroid homomorphism to relate the infinitesimal symmetries of
each structure.

\begin{prop}
If $\JJ$ is a generalized complex structure and $P$ its associated
Poisson structure, then the maps $E\lra T$ defined by $v\mapsto
\pi(v)$ and $v\mapsto\pi(\JJ v)$ both induce homomorphisms
\[
\xymatrix{0\ar[r]&\mathbf{ham}(\JJ)\ar[r]\ar[d]&\mathbf{sym}(\JJ)\ar[r]\ar[d]&H^1(M,L)\ar[d]\ar[r]&0\\
0\ar[r]&\mathbf{ham}(P)\ar[r]&\mathbf{sym}(P)\ar[r]&H^1(M,\Gamma_P)\ar[r]&0}.
\]
from the infinitesimal symmetries of $\JJ$ to the infinitesimal
symmetries of $P$.
\end{prop}
\begin{proof}
Identifying $\textbf{sym}(\JJ)= \ker d_L\cap C^\infty(L^*)$, we see
from Proposition~\ref{morphismo} that $a^*:L^*\lra
\Gamma_P\otimes\CC$ is a morphism of differential complexes.
Identifying $\Gamma_P^*\cong T$, and taking real and imaginary
parts, we obtain morphisms $v\mapsto \pi(v)$, $v\mapsto \pi(\JJ v)$
as required.
\end{proof}

\subsection{Interpolation}

As we saw in section~\ref{typecan}, symplectic structures have type
$0$ while complex structures have type $n$ on a manifold of real
dimension $2n$.  Hence complex and symplectic structures have the
same parity in real dimension $4k$.  We now show that it is possible
to interpolate smoothly between a complex structure and a symplectic
structure through integrable generalized complex structures when $M$
is hyperk\"ahler (or, more generally, holomorphic symplectic).  This
example is also described in~\cite{Hitchin}, using spinors.

Let $M$ be a real manifold of dimension $4k$ with complex structure
$I$ and holomorphic symplectic structure $\sigma = \omega_J +
i\omega_K$, so that $\sigma$ is a nondegenerate closed $(2,0)$-form.
Since $\omega_J$ is of type $(2,0)+(0,2)$, we have $\omega_J I =
I^*\omega_J$, and hence
\begin{equation*}
\begin{pmatrix}
   & -\omega_J^{-1} \\
  \omega_J &  \\
\end{pmatrix}
\begin{pmatrix}
 -I  &  \\
   & I^* \\
\end{pmatrix}
=-\begin{pmatrix}
 -I  &  \\
   & I^* \\
\end{pmatrix}
\begin{pmatrix}
   & -\omega_J^{-1} \\
  \omega_J &  \\
\end{pmatrix},
\end{equation*}
that is, the generalized complex structures $\JJ_{\omega_J}$ and
$\JJ_{I}$ \emph{anticommute}. Hence we may form the one-parameter
family of generalized almost complex structures
\begin{equation*}
\JJ_t=(\sin t) \JJ_I + (\cos t) \JJ_{\omega_J},\ \
t\in[0,\tfrac{\pi}{2}].
\end{equation*}
Clearly $\JJ_t$ is a generalized almost complex structure; we now
check that it is integrable.
\begin{prop}
Let $M$ be a holomorphic symplectic manifold as above.  Then the
generalized almost complex structure $\JJ_t=(\sin t) \JJ_I + (\cos
t) \JJ_{\omega_J}$ is integrable $\forall t\in[0,\tfrac{\pi}{2}]$.
Therefore it is a family of generalized complex structures
interpolating between a symplectic structure and a complex
structure.
\end{prop}
\begin{proof}
Let $B=(\tan t)\omega_K$, a closed 2-form which is well defined
$\forall\ t\in[0,\tfrac{\pi}{2})$.  Noting that $\omega_K I =
I^*\omega_K=\omega_J$, we obtain the following expression:
\begin{equation*}
e^B\JJ_te^{-B}=\left(\begin{matrix}0&-((\sec t)\omega_J)^{-1}\\
(\sec t) \omega_J&0\end{matrix}\right).
\end{equation*}
We conclude from this that for all $t\in[0,\tfrac{\pi}{2})$, $\JJ_t$
is a B-field transform of the symplectic structure determined by
$(\sec t) \omega_J$, and is therefore integrable as a generalized
complex structure; at $t=\tfrac{\pi}{2}$, $\JJ_t$ is purely complex,
and is integrable by assumption, completing the proof.
\end{proof}

\section{Local structure: the generalized Darboux
theorem}\label{darboux}

The Newlander-Nirenberg theorem informs us that an integrable
complex structure on a $2n$-manifold is locally equivalent, via a
diffeomorphism, to $\CC^n$.  Similarly, the Darboux theorem states
that a symplectic structure on a $2n$-manifold is locally
equivalent, via a diffeomorphism, to the standard symplectic
structure $(\RR^{2n},\omega_0)$, where in coordinates
$(x_1,\ldots,x_n, p_1,\ldots,p_n)$,
\begin{equation*}
\omega_0=dx_1\wedge dp_1+\cdots+dx_{n}\wedge dp_n.
\end{equation*}
In this section we prove an analogous theorem for generalized
complex manifolds, describing a local normal form for a
\emph{regular neighbourhood} of a generalized complex manifold.
\begin{defn}
A point $p\in M$ in a generalized complex manifold is called
\emph{regular} when the Poisson structure $P$ is regular at $p$,
i.e. $\type(\JJ)$ is locally constant at $p$.
\end{defn}
By Corollary~\ref{foll}, a generalized complex structure defines, in
a regular neighbourhood $U$, a foliation $\mathcal{F}$ by symplectic
leaves of codimension $2k=2\,\type(\JJ)$, integrating the
distribution $\Delta=\pi(\JJ T^*)$. The complex structure transverse
to $\Delta$ described in Section~\ref{typecan} defines an integrable
complex structure on the leaf space $U/\mathcal{F}$ as we now
describe.
\begin{prop}\label{transcx}
The leaf space $U/\mathcal{F}$ of a regular neighbourhood of a
generalized complex manifold inherits a canonical complex structure.
\end{prop}
\begin{proof}
Let $L\subset E$ be the $+i$-eigenbundle of $\JJ$ and let
$D=\pi_{T\otimes\CC}(L)$ be its projection to the complex tangent
bundle, which is smooth in a regular neighbourhood since $\type(\JJ)
= \dim L\cap (T^*\otimes\CC)$.  Then $E\otimes\CC = L\oplus\overline
L$ implies that $T\otimes\CC = D+\overline D$, while $D\cap
\overline D=\Delta\otimes\CC$.  Since the projection $\pi$ is
bracket-preserving, we see that $D$ is an integrable distribution,
hence $[\Delta, D]\subset D$.  This implies that $D$ descends to an
integrable subbundle $D'\subset T(U/\mathcal{F})\otimes\CC$
satisfying $D'\cap\overline D'=\{0\}$, hence defining an integrable
complex structure on $U/\mathcal{F}$, as required.  This coincides
with the complex structure induced by $\JJ$ on $E/(T^*+\JJ
T^*)=T/\Delta$.
\end{proof}
We now prove that near a regular point, the symplectic structure on
the leaves, together with the complex structure on the leaf space,
completely characterize the generalized complex structure.
\begin{theorem}[Generalized Darboux theorem]\label{thmdar}
Any regular point of type $k$ in a generalized complex manifold has
a neighbourhood which is equivalent, via a diffeomorphism and a
choice of splitting of the Courant algebroid $E$, to the product of
an open set in $\CC^k$ with an open set in the standard symplectic
space $(\RR^{2n-2k},\omega_0)$.
\end{theorem}
\begin{proof}
First choose a local isotropic splitting for the Courant algebroid
so that it is isomorphic, within the neighbourhood $U$, to $(T\oplus
T^*,[\cdot,\cdot]_0)$; this is always possible as long as
$H^3(U,\RR)=0$.

Proposition~\ref{transcx} then guarantees the existence of
holomorphic coordinates $(z_1,\ldots,z_k)$ transverse to the
symplectic foliation in the regular neighbourhood; then a local
generator for the canonical bundle may be chosen, as
in~\eqref{genr}, to be
\begin{equation*}
\rho=e^{B+i\omega}\Omega,
\end{equation*}
where $\Omega=dz_1\wedge\cdots\wedge dz_{k}$ and $B,\omega$ are real
2-forms such that
\begin{equation*}
\omega^{n-k}\wedge\Omega\wedge\overline\Omega \neq 0.
\end{equation*}
Integrability then implies, via Proposition~\ref{locclos}, that
\begin{equation}\label{intcond}
d\rho = e^{B+i\omega}d(B+i\omega)\wedge\Omega = 0.
\end{equation}
The symplectic form $\omega|_\Delta$ along the leaves derives from
the Poisson structure $P$, and hence by Weinstein's normal form for
regular Poisson structures~\cite{WeinPo}, we can find a
leaf-preserving local diffeomorphism
$\varphi:\RR^{2n-2k}\times\CC^k\lra U$ such that
\begin{equation*}
\varphi^*\omega\big|_{\RR^{2n-2k}\times\{pt\}}=\omega_0=dx_1\wedge
dp_1+\cdots+dx_{n-k}\wedge dp_{n-k}.
\end{equation*}

For convenience, let $K=\RR^{2n-2k}$ and $N=\CC^k$, so that
differential forms now have tri-degree $(p,q,r)$ for components in
$\wedge^pK^*\otimes\wedge^qN^*_{1,0}\otimes\wedge^rN^*_{0,1}$.
Furthermore, the exterior derivative decomposes into a sum of three
operators
\begin{equation*}
d=d_\Delta+\del+\delbar,
\end{equation*}
each of degree 1 in the respective component of the tri-grading.
While $\Omega$ is purely of type $(0,k,0)$, the complex 2-form
$A=\varphi^*B+i\varphi^*\omega$ decomposes into six components:
\begin{equation*}\begin{array}{ccc}
 A^{200}  &  &  \\
  A^{110} & A^{101} &  \\
  A^{020} & A^{011} & A^{002}
\end{array}\end{equation*}
Note that only the components $A^{200},A^{101},A^{002}$ act
nontrivially on $\Omega$ in the expression $e^A\Omega$.  Hence we
are free to modify the other three components at will.  Also, note
that the imaginary part of $A^{200}$ is simply $\omega_0$, so that
$d(A^{200}-\overline{A^{200}})=0$, since $\omega_0$ is in constant
Darboux form.

From~\eqref{intcond}, we have $d(B+i\omega)\wedge\Omega=0$, giving
the following four equations:
\begin{align}
\delbar A^{002}&=0\label{first}\\
\delbar A^{101}+d_\Delta A^{002}&=0\label{second}\\
\delbar A^{200}+d_\Delta A^{101}&=0\label{third}\\
d_{\Delta} A^{200}&=0.\label{fourth}
\end{align}

We will now endeavour to modify $A$ so that
$\varphi^*\rho=e^A\Omega$ is unchanged but $A$ is replaced with
$\tilde{A}=\tilde{B}+\tfrac{1}{2}(A^{200}-\overline{A^{200}})$,
where $\tilde{B}$ is a real closed 2-form.  This would demonstrate
that
\begin{equation*}
\varphi^*\rho=e^{\tilde B + i\omega_0}\Omega,
\end{equation*}
i.e. $\rho$ is equivalent, via a diffeomorphism and $B$-field
symmetry, to the product of a symplectic with a complex structure.
The $B$-field transform is simply a change in the original splitting
for $(T\oplus T^*,[\cdot,\cdot]_0)$.

In order to preserve $\varphi^*\rho$, the most general form for
$\tilde B$ is
\begin{equation*}
\tilde B =
\tfrac{1}{2}(A^{200}+\overline{A^{200}})+A^{101}+\overline{A^{101}}+A^{002}+\overline{A^{002}}
+ C,
\end{equation*}
where $C$ is a real 2-form of type $(011)$.  Then clearly
$\varphi^*\rho = e^{\tilde{B}+i\omega_0}\Omega$.  Requiring that
$d\tilde B=0$ imposes two constraint equations:
\begin{align}
(d\tilde B)^{012}&=\del A^{002}+\delbar C = 0.\label{cond1}\\
(d\tilde B)^{111}&=\del A^{101}+\overline{\del A^{101}} + d_\Delta
C=0\label{cond2}
\end{align}
The question then becomes whether we can find a real $(011)$-form
$C$ such that these equations are satisfied.  The following are all
local arguments, making repeated use of the Dolbeault lemma.
\begin{itemize}
\item From equation (\ref{first}) we obtain that
$A^{002}=\delbar\alpha$ for some $(001)$-form $\alpha$.  Then
condition (\ref{cond1}) is equivalent to $\delbar(C-\del\alpha)=0$,
whose general solution is
\[C=\del\alpha+\overline{\del\alpha} + i\del\delbar\chi\]for any real
function $\chi$. We must now check that it is possible to choose
$\chi$ so that the second condition (\ref{cond2}) is satisfied by
this $C$.

\item From equation (\ref{second}) we obtain that
$\delbar(A^{101}-d_\Delta\alpha)=0$, implying that
$A^{101}=d_\Delta\alpha+\delbar\beta$ for some $(100)$-form $\beta$.
Condition (\ref{cond2}) then is equivalent to the fact that
\begin{equation*}
-id_\Delta\del\delbar\chi=\del\delbar(\beta-\overline\beta),
\end{equation*}
which can be solved (for the unknown $\chi$) if and only if the
right hand side is $d_\Delta$-closed.  From equation (\ref{third})
we see that $\delbar(A^{200}-d_\Delta\beta)=0$, showing that
$A^{200}=d_\Delta\beta + \delta$, where $\delta$ is a
$\delbar$-closed $(200)$-form. Hence
\begin{equation*}
d_\Delta\del\delbar(\beta-\overline\beta)=\del\delbar(A^{200}-\overline{A^{200}}),
\end{equation*}
and the right hand side vanishes precisely because
$A^{200}-\overline{A^{200}}=2\omega_0$, which is closed.  Hence
$\chi$ may be chosen to satisfy condition (\ref{cond2}), and so we
obtain a closed 2-form $\tilde B$.
\end{itemize}
\end{proof}

\subsection{Type jumping}\label{examples}

While Theorem~\ref{thmdar} fully characterizes generalized complex
structures in regular neighbourhoods, it remains an essential
feature of the geometry that the \emph{type} of the structure may
vary throughout the manifold.  The most generic type is zero, when
there are only symplectic directions and the Poisson structure $P$
has maximal rank.  The type may jump up along closed subsets, has
maximal value $n=\tfrac{1}{2}\dim_\RR M$, and has fixed parity
throughout the manifold.  We now present a simple example of a
generalized complex structure on $\RR^4$ which is of symplectic type
$(k=0)$ outside a codimension 2 hypersurface and jumps up to complex
type $(k=2)$ along the hypersurface.

Consider the differential form
\begin{equation}\label{tc}
\rho=z_1 + dz_1\wedge dz_2,
\end{equation}
where $z_1,z_2$ are the standard coordinates on $\CC^2\cong\RR^4$.
Along $z_1=0$, we have $\rho=dz_1\wedge dz_2$ and so it generates
the pure spinor line  corresponding to the standard complex
structure. Whenever $z_1\neq 0$, $\rho$ may be rewritten as follows:
\begin{equation*}
\rho=z_1 e^{\frac{dz_1\wedge dz_2}{z_1}}.
\end{equation*}
Therefore away from $z_1=0$, $\rho$ generates the canonical line
bundle of the $B$-field transform of the symplectic form $\omega$,
where
\begin{equation*}
B+i\omega=z_1^{-1}dz_1\wedge dz_2.
\end{equation*}
Hence, algebraically the form $\rho$ defines a generalized almost
complex structure which is generically of type $0$ but jumps to type
$2$ along $z_1=0$.

To verify the integrability of this structure, we take the exterior
derivative:
\begin{equation*}
d\rho=dz_1=i_{-\partial_{z_2}}(z_1+dz_1\wedge dz_2)=
(-\partial_{z_2})\cdot\rho,
\end{equation*}
showing that $\rho$ indeed satisfies the integrability condition of
Theorem~\ref{deldelbar}, and defines a generalized complex structure
on all of $\RR^4$.  In this case it is easy to see that although the
canonical line bundle is topologically trivial, it does not admit a
closed, nowhere-vanishing section. Hence the generalized complex
structure is not generalized Calabi-Yau.

In the next chapter we will produce more general examples of the
jumping phenomenon, and on compact manifolds as well.  However we
indicate here that the simple example above was used in~\cite{GMJDG}
to produce, via a surgery on a symplectic 4-manifold, an example of
a compact, simply-connected generalized complex 4-manifold which
admits neither complex nor symplectic structures.

\section{Deformation theory}\label{deform}

In the deformation theory of complex manifolds developed by Kodaira,
Spencer, and Kuranishi, one begins with a compact complex manifold
$(M,J)$ with holomorphic tangent bundle $\TT$, and constructs an
analytic subvariety $\zed\subset H^1(M,\TT)$ (containing $0$) which
is the base space of a family of deformations $\MM=\{\eps(z)\ :\
z\in \zed,\ \eps(0)=0\}$ of the original complex structure $J$.
This family is \emph{locally complete} (also called
\emph{miniversal}), in the sense that any family of deformations of
$J$ can be obtained, up to equivalence, by pulling $\MM$ back by a
map $f$ to $\zed$, as long as the family is restricted to a
sufficiently small open set in its base.

The subvariety $\zed\subset H^1(M,\TT)$ is defined as the zero set
of a holomorphic map $\Phi:H^1(M,\TT)\rightarrow H^2(M,\TT)$, and so
the base of the miniversal family is certainly smooth when this
obstruction map vanishes.

In this section we extend these results to the generalized complex
setting, following the method of Kuranishi~\cite{Kuranishi}.  In
particular, we construct, for any generalized complex manifold, a
locally complete family of deformations.  We then proceed to produce
new examples of generalized complex structures by deforming known
ones.

\subsection{Lie bialgebroids and the deformation complex}

The generalized complex structure $\JJ$ on the exact Courant
algebroid $E$ is determined by its $+i$-eigenbundle $L\subset
E\otimes\CC$ which is isotropic, satisfies
$L\cap\overline{L}=\{0\}$, and is Courant involutive. Recall that
since $E\otimes\CC=L\oplus\overline L$, we use the natural metric
$\IP{,}$ to identify $\overline L$ with $L^*$.

To deform $\JJ$ we will vary $L$ in the Grassmannian of maximal
isotropics. Any maximal isotropic having zero intersection with
$\overline{L}$ (this is an open set containing $L$) can be uniquely
described as the graph of a homomorphism $\epsilon:L\lra \overline
L$ satisfying $\IP{\epsilon X,Y }+\IP{X,\epsilon Y}=0\ \ \forall
X,Y\in C^\infty(L)$, or equivalently $\epsilon\in C^\infty(\wedge^2
L^*)$. Therefore the new isotropic is given by
\begin{equation*}
L_\eps=(1+\eps)L=\{u+i_u\eps\ :\ u\in L\}.
\end{equation*}
As the deformed $\JJ$ is to remain real, we must have $\overline
L_\epsilon=(1+\overline\epsilon)\overline L$.  Now we observe that
$L_\epsilon$ has zero intersection with its conjugate if and only if
the endomorphism we have described on $L\oplus L^*$, namely
\begin{equation}\label{matrixa}
A_\epsilon=\left(\begin{matrix}1&\overline\epsilon\\\epsilon&1\end{matrix}\right),
\end{equation}
is invertible; this is the case for $\epsilon$ in an open set around
zero.

So, providing $\epsilon$ is small enough,
$\JJ_\epsilon=A_\epsilon\JJ A_\epsilon^{-1}$ is a new generalized
almost complex structure, and all nearby almost structures are
obtained in this way.  Note that while $A_\epsilon$ itself is not an
orthogonal transformation, of course $\JJ_\epsilon$ is.

To describe the condition on $\epsilon\in C^\infty(\wedge^2 L^*)$
which guarantees that $\JJ_\epsilon$ is integrable, we observe the
following.  Since $L^*=\overline L$, we have not only an elliptic
differential complex (by Proposition~\ref{ellipcx})
\[
\left( C^\infty(\wedge^\bullet L^*), d_L \right),
\]
but also a Lie algebroid structure on $L^*$ coming from the Courant
bracket on $\overline L$. In fact, by a theorem of
Liu-Weinstein-Xu~\cite{LWX}, the differential is a derivation of the
bracket and we obtain the structure of a Lie bialgebroid in the
sense of Mackenzie-Xu~\cite{MackXu}, also known as a differential
Gerstenhaber algebra.

\begin{theorem}[\cite{LWX}, Theorem 2.6]\label{bial}
Let $E$ be an exact Courant algebroid and $E=L\oplus L'$ for Dirac
structures $L,L'$. Then $L'=L^*$ using the inner product, and the
dual pair of Lie algebroids $(L,L^*)$ defines a Lie bialgebroid,
i.e.
\begin{equation*}
d_L[a,b]=[d_L a,b]+[a,d_L b],
\end{equation*}
for $a,b\in C^\infty(L^*)$, where $[\cdot,\cdot]$ is extended in the
Schouten sense to $C^\infty(\wedge^\bullet L^*)$.  Therefore the
data
\[
\left(C^\infty(\wedge^\bullet L^*), d_L, [\cdot,\cdot]\right)
\]
define a differential Gerstenhaber algebra.
\end{theorem}

Interpolating between the examples~\ref{symp} and~\ref{pois},
Liu-Weinstein-Xu~\cite{LWX} also prove that, under the assumptions
of the previous theorem, the graph $L_\eps$ of a section $\eps\in
C^\infty(\wedge^2 L^*)$ defines an integrable Dirac structure if and
only if it satisfies the Maurer-Cartan equation.
\begin{theorem}[\cite{LWX}, Theorem 6.1]\label{masterequation}
The almost Dirac structure $L_\eps$, for $\eps\in
C^\infty(\wedge^2L^*)$, is integrable if and only if $\eps$
satisfies the Maurer-Cartan equation
\begin{equation}\label{MC}
d_L\eps+\tfrac{1}{2}[\eps,\eps]=0.
\end{equation}
Here $d_L:C^\infty(\wedge^k L^*)\rightarrow
C^\infty(\wedge^{k+1}L^*)$ and $[\cdot,\cdot]$ is the Lie algebroid
bracket on $L^*$.
\end{theorem}

Therefore we conclude that the deformed generalized almost complex
structure $\JJ_\epsilon$ is integrable if and only if $\epsilon$
satisfies the Maurer-Cartan equation~\eqref{MC}.  We may finally
define a smooth family of deformations of the generalized complex
structure $\JJ$.  We are only interested in ``small'' deformations.
\begin{defn}
Let $U$ be an open disk containing the origin of a
finite-dimensional vector space. A \emph{smooth family of
deformations} of $\JJ$ over $U$ is a family of sections $\eps(u)\in
C^\infty(\wedge^2 L^*)$, smoothly varying in $u\in U$, with
$\eps(0)=0$, such that~\eqref{matrixa} is invertible, and satisfying
the Maurer-Cartan equation~\eqref{MC} for each $u\in U$. Two such
families $\eps_1(u),\eps_2(u)$ are \emph{equivalent} if
$F_u(L_{\eps_1}(u))=L_{\eps_2}(u)$ for all $u\in U$, where $F_u$ is
a smooth family of Courant automorphisms with $F_0=\id$.
\end{defn}

The space of solutions to~\eqref{MC} is infinite-dimensional,
however due to the action of the group of Courant automorphisms we
are able, as in the case of complex manifolds, to take a suitable
quotient, forming a finite-dimensional locally complete family. To
obtain this finite-dimensional moduli space of deformations, it will
suffice to consider equivalences $F_u$ which are families of exact
Courant automorphisms in the sense of Definition~\ref{exaut},
generated by \emph{time-independent} derivations $\ad (v(u))$ given
by a smooth family of sections $v(u)\in C^\infty(E)$, $u\in U$.  A
similar situation occurs in the case of deformations of complex
structure.

Suppose $v\in C^\infty(E)$ and let $F_v^1$ denote its time-1 flow
defined by~\eqref{oneflow}, so that in a splitting for $E$ with
curvature $H$, we have $v = X+\xi$ and by~\eqref{intbee},
\begin{equation}\label{einflow}
F^1_v= \varphi^1_* e^{B_1},\ \ \ B_1=\int_{0}^1 \varphi^*_s(i_XH
+d\xi)\ ds,
\end{equation}
where $\varphi^t_*$ is the flow of the vector field $X$. The Courant
isomorphism $F^1_v$ acts on generalized complex structures, taking a
given deformation $L_\eps$ to $F^1_v(L_\eps)$. If $v$ has
sufficiently small 1-jet, then $F^1_v(L_\eps)$ may be expressed as
$L_{\eps'}$ for another section $\eps'\in C^\infty(\wedge^2 L^*)$,
and we denote it
 $F^1_v(\eps):=\eps'$. We now determine an approximate formula for $F^1_v(\eps)$ in
terms of $(\eps,v)$.

\begin{prop}
Let $\JJ$ be a generalized complex structure with $+i$-eigenbundle
$L\subset E\otimes\CC$, and let $\eps\in C^\infty(\wedge^2 L^*)$ be
such that~\eqref{matrixa} is invertible. Then for $v\in C^\infty(E)$
with sufficiently small 1-jet, the time-1 flow~\eqref{einflow}
satisfies
\begin{equation}\label{transdefm}
F^1_v(\eps)=\eps + d_Lv^{0,1}+R(\eps,v),
\end{equation}
where $v = v^{1,0} +v^{0,1}$ according to the splitting $L\oplus
L^*$, and $R$ satisfies
\begin{equation*}
R(t\eps,tv)=t^2\widetilde R(\eps,v,t),
\end{equation*}
where $\widetilde R(\eps,v,t)$ is smooth in $t$ for small $t$.
\end{prop}
\begin{proof}
Define $\eps(s,t)$ for $s,t\in \RR$ by
\begin{equation}\label{plan}
\eps(s,t) = F^1_{tv}(s\eps),
\end{equation}
so that $\eps=\eps(1,1)$ and $\eps(0)=0$.  We first compute the
derivatives of~\eqref{plan} at $s=t=0$.  The derivative in $s$ is
easily computed:
\begin{equation*}
\left.\frac{\del \eps(s,t)}{\del s}\right|_{(0,0)} = \frac{\del
(s\eps)}{\del s} = \eps.
\end{equation*}
The derivative in $t$ may be computed using the property of flows
that $F^1_{tv} = F^t_{v}$, together with the fact that the flow
$F^t_{v}$ is generated by the adjoint action $\ad(v)=[v,\cdot]$ of
$v$ on the Courant algebroid.  Using the properties of the Courant
bracket, we obtain, for $y,z\in C^\infty(L)$,
\begin{equation*}
\left.\frac{\del \eps(s,t)}{\del t}\right|_{(0,0)}(y,z) = \frac{\del
F^t_v(0)}{\del t}(y,z) = \IP{-[v,y],z} = d_Lv^{0,1}(y,z).
\end{equation*}
By Taylor's theorem we obtain
\[
F^1_{tv}(s\eps) = s\eps + td_L v^{0,1} + r(s,t,\eps,v),
\]
where $r$ is smooth of order $O(s^2,st,t^2)$ at zero.  Setting
$R(\eps,v)=r(1,1,\eps,v)$, we obtain the result, since clearly
$r(1,1,t\eps,tv)=r(t,t,\eps,v)$ is of order $O(t^2)$.
\end{proof}

Whereas the Maurer-Cartan equation~\eqref{MC} indicates that,
infinitesimally, deformations of generalized complex structure lie
in  $\ker d_L\subset C^\infty(\wedge^2 L^*)$, the previous
Proposition shows us that, infinitesimally, deformations which
differ by sections which lie in the image of $d_L$ are equivalent.
Hence we expect the tangent space to the moduli space to lie in the
Lie algebroid cohomology $H^2(M,L)$, which by ellipticity is
finite-dimensional for $M$ compact. We now develop the Hodge theory
required to prove this assertion.

We follow the usual treatment of Hodge theory as described
in~\cite{Wells}. Choose a Hermitian metric on the complex Lie
algebroid $L$ and let $|\varphi|_k$ be the $L^ 2_k$ Sobolev norm on
sections $\varphi\in C^\infty(\wedge^p L^*)$ induced by the metric.
We then have the elliptic, self-adjoint Laplacian
\begin{equation*}
\Delta_L=d_Ld^*_L+d^*_Ld_L.
\end{equation*}
Let $\Har^p$ be the space of $\Delta_L$-harmonic forms, which is
isomorphic to $H^p(M,L)$ by the standard argument, and let $H$ be
the orthogonal projection of $C^\infty(\wedge^p L)$ onto the closed
subspace $\Har^p$.  Also, let $G$ be the Green smoothing operator
quasi-inverse to $\Delta_L$, i.e. $G\Delta + H = \mathrm{Id}$ and
\begin{equation*}
G:L^2_k\rightarrow L^2_{k+2}.
\end{equation*}
We will find it useful, as Kuranishi did, to define the
once-smoothing operator
\begin{equation*}
Q=d_L^*G:L^2_k\rightarrow L^2_{k+1},
\end{equation*}
which then satisfies
\begin{align}\label{Qrel}
\mathrm{Id}&=H+d_LQ+Q d_L,\\
Q^2=d^*_LQ&=Qd^*_L=HQ=QH=0.\nonumber
\end{align}
We now have the algebraic and analytical tools required to prove a
direct analog of Kuranishi's theorem for generalized complex
manifolds.

\subsection{The deformation theorem}

\begin{theorem} Let $(M,\JJ)$ be a compact generalized complex manifold.  There exists an open neighbourhood $U\subset H^2(M,L)$
containing zero, a smooth family $\widetilde\MM=\{\eps(u)\ :\ u\in
U,\ \eps(0)=0\}$ of generalized almost complex deformations of
$\JJ$, and an analytic obstruction map $\Phi:U\rightarrow H^3(M,L)$
with $\Phi(0)=0$ and $d\Phi(0)=0$, such that the deformations in the
sub-family $\MM=\{\eps(z)\ :\ z\in\zed=\Phi^{-1}(0)\}$  are
precisely the integrable ones. Furthermore, any sufficiently small
deformation $\eps$ of $\JJ$ is equivalent to at least one member of
the family $\MM$. In the case that the obstruction map vanishes,
$\MM$ is a smooth locally complete family.
\end{theorem}
\begin{proof} The proof is divided into two parts: first, we
construct a smooth family $\widetilde\MM$, and show it contains the
family of integrable deformations $\MM$ defined by the map $\Phi$;
second, we describe its miniversality property.  We follow the paper
of Kuranishi~\cite{Kuranishi} closely, where more details can be
found.

{\bf Part I:} For sufficiently large $k$, $L^2_k(M,\RR)$ is a Banach
algebra (see~\cite{Palais}), and the map $f:\eps\mapsto
\eps+\tfrac{1}{2}Q[\eps,\eps]$ extends to a smooth map
\begin{equation*}
f:L^2_k(\wedge^2 L^*)\lra L^2_k(\wedge^2 L^*),
\end{equation*}
whose derivative at the origin is the identity mapping. By the
inverse function theorem, $f^{-1}$ maps a neighbourhood of the
origin in $L^2_k(\wedge^2 L^*)$ smoothly and bijectively to another
neighbourhood of the origin. Hence, for sufficiently small
$\delta>0$, the finite-dimensional subset of harmonic sections,
\begin{equation*}
U=\{u\in\Har^2<L^2_k(\wedge^2L^*)\ : \ |u|_k<\delta\},
\end{equation*}
defines a family of sections as follows:
\begin{equation*}
\widetilde M=\{\eps(u)=f^{-1}(u)\ :\ u\in U\},
\end{equation*}
where $\eps(u)$ depends smoothly (in fact, holomorphically) on $u$,
and satisfies $f(\eps(u))=u$. Applying the Laplacian to this
equation, we obtain
\begin{equation*}
\Delta_L\eps(u) + \tfrac{1}{2}d^*_L[\eps(u),\eps(u)]=0.
\end{equation*}
This is a quasi-linear elliptic PDE, and by a result of
Morrey~\cite{Morrey}, we conclude that the solutions $\eps(u)$ of
this equation are actually smooth, i.e.
\begin{equation*}
\eps(u)\in C^\infty(\wedge^2 L^*).
\end{equation*}
Hence we have constructed a smooth family of generalized almost
complex deformations of $\JJ$, over an open set $U\subset\Har^2\cong
H^2(M,L)$.

We now ask which of these deformations satisfy the Maurer-Cartan
equation~\eqref{MC}.  By definition of $\eps(u)$, and using
~\eqref{Qrel}, we obtain
\begin{align*}
d_L\eps(u)+\tfrac{1}{2}[\eps(u),\eps(u)]&=-\tfrac{1}{2}d_LQ[\eps(u),\eps(u)]+\tfrac{1}{2}[\eps(u),\eps(u)]\\
&=\tfrac{1}{2}(Qd_L+H)[\eps(u),\eps(u)].
\end{align*}
Since the images of $Q$ and $H$ are $L^2$-orthogonal, we see that
$\eps(u)$ is integrable if and only if
$H[\eps(u),\eps(u)]=Qd_L[\eps(u),\eps(u)]=0$.  We now refer to the
argument of Kuranishi~\cite{Kuranishi} which, using the
compatibility between $[\cdot,\cdot]$ and $d_L$, shows that
$Qd_L[\eps(u),\eps(u)]$ vanishes when $H[\eps(u),\eps(u)]$ does.

Hence, $\eps(u)$ is integrable precisely when $u$ lies in the
vanishing set of the analytic mapping $\Phi:U\rightarrow H^3_L(M)$
defined by
\begin{equation}\label{obstruc}
\Phi(u)=H[\eps(u),\eps(u)].
\end{equation}
Note furthermore that $\Phi(0)=d\Phi(0)=0$.

{\bf Part II:} For the second part of the proof, we give an
alternative characterisation of the family $\MM=\{\eps(z)\ :\ z\in
\zed=\Phi^{-1}(0) \}$.  We claim that $\MM$ is actually a
neighbourhood around zero in the set
\begin{equation*}
\MM'=\left\{ \eps\in C^\infty(\wedge^2 L^*)\ :\
d_L\eps+\tfrac{1}{2}[\eps,\eps]=0,\ \ d^*_L\eps=0\right\}.
\end{equation*}
To show this, let $\eps(u)\in \MM$. Then since
$\eps(u)=u-\tfrac{1}{2}Q[\eps(u),\eps(u)]$ and $d^*_LQ=0$, we see
that $d^*_L\eps(u)=0$, showing that $\MM\subset\MM'$.  Conversely,
let $\eps\in\MM'$.  Then since $d^*_L\eps=0$, applying $d^*_L$ to
the equation $d_L\eps+\tfrac{1}{2}[\eps,\eps]=0$ we obtain
$\Delta_L\eps+\tfrac{1}{2}d^*_L[\eps,\eps]=0$, and applying Green's
operator we see that $\eps+\tfrac{1}{2}Q[\eps,\eps]=H\eps$, i.e.
$F(\eps)=H\eps\in\Har^2$, proving that a small open set in $\MM'$ is
contained in $\MM$, completing the argument.

We now show that every sufficiently small deformation of the
generalized complex structure is equivalent to one in our
finite-dimensional family $\MM$.  Let $P\subset C^\infty(L^*)$ be
the $L^2$ orthogonal complement of $\ker d_L\subset C^\infty(L^*)$,
or in other words, sections in the image of $d^*_L$. We show that
there exist neighbourhoods of the origin $V\subset C^\infty(\wedge^2
L^*)$ and $W\subset P$ such that for any $\eps\in V$ there is a
unique $v\in C^\infty(E)$ such that $v^{0,1}\in W$ and the time-1
flow $F^1_v(\eps)$ satisfies
\begin{equation}\label{flowvan}
d^*_LF^1_{v}(\eps)=0.
\end{equation}
This would imply that any sufficiently small solution to
$d_L\eps+\tfrac{1}{2}[\eps,\eps]=0$ is equivalent to another
solution $\eps'$ such that $d^*_L\eps'=0$, i.e. a solution in
$\MM'=\MM$. Extended to smooth families, this result would prove
local completeness.

We see from~\eqref{transdefm} that~\eqref{flowvan} holds if and only
if
\begin{equation*}
d^*_L\eps +d^*_Ld_Lv^{0,1}+d^*_LR(\eps,v)=0.
\end{equation*}
Assuming $v^{0,1}\in P$, we see that $d^*_Lv^{0,1}=Hv^{0,1}=0$, so
that
\begin{equation}\label{qell}
d^*_L\eps +\Delta_Lv^{0,1}+d^*_LR(\eps,v)=0.
\end{equation}
Applying the Green operator $G$, we obtain
\begin{equation}\label{vaneq}
v^{0,1}+Q\eps +QR(\eps,v)=0.
\end{equation}
By the definition~\eqref{transdefm} of $R(\eps,v)$, the map
\begin{equation*}
F:(\eps,v^{0,1})\mapsto v^{0,1}+Q\eps +QR(\eps,v)
\end{equation*}
is continuous from a neighbourhood of the origin $V_0\times W_0$ in
$C^\infty(\wedge^2 L^*)\times P$ to $P$, where all spaces are
endowed with the $L^2_k$ norm, $k$ sufficiently large. Also, the
derivative of $F$ with respect to $v^{0,1}$ (at $0$) is the identity
map.  Therefore by the implicit function theorem, there are
neighbourhoods $V\subset V_0, W_1\subset \widehat{W_0}$ such that
given $\eps\in V$, equation~\eqref{vaneq}, i.e. $F=0$, is satisfied
for a unique $v^{0,1}\in W_1$, and which depends smoothly on
$\eps\in V$. Furthermore, since $\eps\in V$ is itself smooth, the
unique solution $v$ satisfies the quasi-linear elliptic
PDE~\eqref{qell}, implying that $v$ is smooth as well, hence
$v^{0,1}$ lies in the neighbourhood $W=W_1\cap P$. Therefore we have
shown that every sufficiently small deformation of the generalized
complex structure is equivalent to one in our finite-dimensional
family $\MM$.

If the obstruction map $\Phi$ vanishes, so that $\MM$ is a smooth
family, then given any other smooth family $\MM_S=\{\eps(s)\ :\ s\in
S,\ \eps(s_0)=0\}$ with basepoint $s_0\in S$, the above argument
provides, for $s$ in some neighbourhood $T$ of $s_0$, a smooth
family of sections $v(s)\in C^\infty(E)$ whose time-1 flow takes
each $\eps(s)$ to $\eps(f(s))$, $f(s)\in U\subset H^2(M,L)$.  This
defines a smooth map $f:T\rightarrow U$, $f(s_0)=0$, such that
$f^*\MM=\MM_S$. Thus we establish that $\MM$ is a locally complete
family of deformations.
\end{proof}
\begin{remark}
The natural complex structure on $H^2(M,L)$ and on the vanishing set
of the holomorphic obstruction map $\Phi$ raises the question of
whether there is a notion of holomorphic family of generalized
complex structures.  There is: if $S$ is a complex manifold then a
holomorphic family of generalized complex structures on $M$ is a
generalized complex structure on $M\times S$ which can be pushed
down, or \emph{reduced} in the sense of~\cite{BCG}, via the
projection to yield the complex structure on $S$. The family $\MM$
can be shown to define such a holomorphic family, since the
constructed family $\eps(u)$ depends holomorphically on $u$.
\end{remark}

\subsection{Examples of deformed structures}\label{deformex}

Consider deforming a compact complex manifold $(M,J)$ as a
generalized complex manifold.  Recall that the associated Lie
algebroid is $L=T_{0,1}\oplus T^*_{1,0}$, so the deformation complex
is simply the holomorphic multivector Dolbeault complex
\begin{equation*}
\left(\Omega^{0,\bullet}(\wedge^\bullet T_{1,0}),\delbar\right).
\end{equation*}
The base of the Kuranishi family therefore lies in the
finite-dimensional vector space
\begin{equation*}
H^2(M,L)=\oplus_{p+q=2}H^q(M,\wedge^p T_{1,0}),
\end{equation*}
whereas the image of the obstruction map lies in
\begin{equation*}
H^3(M,L)=\oplus_{p+q=3}H^q(M,\wedge^p T_{1,0}).
\end{equation*}
In this way, generalized complex manifolds provide a geometrical
interpretation of the ``extended complex deformation space'' defined
by Kontsevich and Barannikov~\cite{KontsevichBarannikov}. Any
deformation $\eps$ has three components
\begin{equation*}
\beta\in H^0(M,\wedge^2T_{1,0}),\ \ \ \varphi\in H^1(M,T_{1,0}),\ \
\ B\in H^2(M,\mathcal{O}).
\end{equation*}
The component $\varphi$ is a usual deformation of the complex
structure, as described by Kodaira and Spencer. The component $B$
represents a residual action by cohomologically nontrivial $B$-field
transforms; these do not affect the type. The component $\beta$,
however, is a new deformation for complex manifolds. Setting
$B=\varphi=0$, the integrability condition reduces to
\begin{equation*}
\delbar\beta + \tfrac{1}{2}[\beta,\beta]=0,
\end{equation*}
which is satisfied if and only if the bivector $\beta$ is
holomorphic and Poisson.  Writing $\beta = \tfrac{-1}{4}(Q+iP)$ for
$Q,P$ real bivectors of type $(0,2)+(2,0)$ such that $Q=PJ^*$, we
may explicitly determine the deformed generalized complex structure:
\begin{equation*}
\JJ_\beta = e^{\beta +
\overline\beta}\JJ_Je^{-(\beta+\overline\beta)}=\begin{pmatrix}J&P\\
&-J^*\end{pmatrix}.
\end{equation*}
In this way, we obtain a new class of generalized complex manifolds
with type controlled by the rank of the holomorphic Poisson bivector
$\beta$.

\begin{example}[Deformed generalized complex structure on $\CC
P^2$]\label{defmcp2} For $\CC P^2$, $\wedge^2
T_{1,0}=\mathcal{O}(3)$, and for dimensional reasons, any
holomorphic bivector $\beta\in H^0(M,\mathcal{O}(3))$ is
automatically Poisson.  Hence any holomorphic section of
$\mathcal{O}(3)$ defines an integrable deformation of the complex
structure into a generalized complex structure.  Since
$H^1(T_{1,0})=H^2(\mathcal{O})=0$, we may conclude from the
arguments above that the locally complete family of deformations is
smooth and of complex dimension 10.  However one can also check that
the obstruction space vanishes in this case, by the Bott formulae.

The holomorphic Poisson structure $\beta$ has maximal rank outside
its vanishing locus, which must be a cubic curve $C$. Hence the
deformed generalized complex structure is of $B$-symplectic type
(type 0) outside $C$ and of complex type (type 2) along the cubic.
The complexified symplectic form $B+i\omega = \beta^{-1}$ is
singular along $C$.  We therefore have an example of a compact
generalized complex manifold exhibiting type change along a
codimension 2 subvariety.
\end{example}

\begin{example}
One can deform the standard complex structure on $\CC^2$ by the
holomorphic bivector
\begin{equation*}
\beta=z_1\del_{z_1}\wedge\del_{z_2},
\end{equation*}
where $z_1,z_2$ are the usual complex coordinates. Applying a
$\beta$-transform to the spinor $\Omega = dz_1\wedge dz_2$ defining
the usual complex structure, we obtain
\begin{equation*}
e^{\beta}\Omega= dz_1\wedge dz_2 + z_1,
\end{equation*}
which is precisely the example~\eqref{tc}.  We see now that it is
actually a deformation of the usual complex structure by a
holomorphic Poisson structure.
\end{example}

Note that while holomorphic Poisson bivectors may be thought of as
infinitesimal noncommutative deformations in the sense of
quantization of Poisson structures, we are viewing them here as
genuine (finite) deformations of the generalized complex structure.
For more details about this distinction and its consequences,
see~\cite{Kap2},\cite{GMN}.

\section{Generalized complex branes}\label{branes}

In this section we introduce the natural ``sub-objects'' of
generalized complex manifolds, generalizing both holomorphic
submanifolds of a complex manifold and Lagrangian submanifolds in
symplectic geometry.  In fact, even in the case of a usual
symplectic manifold, there are generalized complex branes besides
the Lagrangian ones: we show these are the \emph{coisotropic
A-branes} discovered by Kapustin and Orlov~\cite{Kapustin}.

As has been emphasized by physicists, a geometric description of a
brane in $M$ involves not only a submanifold $\iota: S\lra M$ but
also a vector bundle supported on it; in cases where a nontrivial
$S^1$-gerbe $G$ is present one replaces the vector bundle by an
\emph{object} (``twisted vector bundle'') of the pullback gerbe
$\iota^*G$. Since the Courant bracket captures the differential
geometry of the gerbe, we obtain a convenient description of branes
in terms of generalized geometry. For simplicity we shall restrict
our attention to branes supported on loci where the pullback gerbe
$\iota^*G$ is trivializable.

We begin by phrasing the definition of a gerbe trivialization in
terms of the Courant bracket. Recall that $\iota^*E$ denotes the
pullback of exact Courant algebroids, defined in the Appendix.
\begin{defn}
Let $E$ be an exact Courant algebroid on $M$ and let $\iota:S\lra M$
be a submanifold.  A (Courant) \emph{trivialization} of $E$ along
$S$ consists of a bracket-preserving isotropic splitting $s:TS\lra
\iota^*E$ inducing an isomorphism
\[
s+\pi^*:(TS\oplus T^*S,[\cdot,\cdot]_0)\lra \iota^*E.
\]
\end{defn}

If an isotropic splitting $\tilde s:TM\lra E$ is chosen, with
curvature $H\in\Omega^3_{cl}(M)$, then $\iota^*E$ inherits a
splitting with curvature $\iota^*H$, and any trivialization
$(\iota,s)$ is characterized by the difference $s-\iota^*\tilde
s=F\in \Omega^2(S)$, which satisfies
\begin{equation}\label{trivg}
\iota^*H = dF.
\end{equation}
Therefore the gerbe curvature is exact when pulled back to $S$.
Indeed, we obtain a generalized pullback morphism $\rho\mapsto
e^F\wedge \iota^*\rho$, defining a map from the twisted de Rham
complex of $M$ to the usual de Rham complex of $S$:
\begin{equation*}
\xymatrix{(\Omega^\bullet(M),d_H)\ar[r]^{e^F\iota^*}&
(\Omega^\bullet(S),d)}.
\end{equation*}
This may be viewed as the the image under the Chern character of a
morphism from the twisted K-theory of $M$ to the usual K-theory of
$S$.

To avoid confusion, let $E|_S$ denote the restriction of the bundle
$E$ to $S$, as opposed to $\iota^*E=K^\perp/K$, for $K=\Ann(TS)$,
which defines the pullback Courant algebroid over $S$.  The
trivialization $(\iota,s)$ defines a maximal isotropic subbundle
$s(TS)\subset \iota^*E$.  Further, the quotient map $q:K^\perp\lra
K^\perp/K$ determines a bijection taking maximal isotropic
subbundles $L\subset \iota^*E$ to maximal isotropic subbundles
$q^{-1}(L)\subset E|_S$ contained in $K^\perp$.
\begin{defn}
The \emph{generalized tangent bundle}  to the trivialization
$\mathcal{L}=(\iota,s)$ of $E$ is the maximal isotropic subbundle
$\tau_\LL\subset E|_S$ defined by $\tau_\LL=q^{-1}(s(TS))$.
\end{defn}
Note that $\Ann(TS)=N^*S$, so that $\tau_\LL$ is actually an
extension of the tangent bundle by the conormal bundle:
\begin{equation}\label{exttau}
\xymatrix{0\ar[r]&N^*S\ar[r]& \tau_\LL\ar[r]^\pi& TS\ar[r]&0}
\end{equation}
If a splitting $\tilde s:TM\lra E$ is chosen, with $s-\iota^*\tilde
s = F\in \Omega^2(S)$ as in~\eqref{trivg}, then $\tau_\LL$ has the
explicit form
\begin{equation}\label{exptau}
\tau_\LL = \{X+\eta\in TS\oplus T^*M\ :\ \iota^*\eta=i_XF\}.
\end{equation}
Comparing this with Proposition~\ref{regular}, we obtain the
following canonical example of a Courant trivialization:
\begin{example}\label{extbra}
Let $L\subset (TM\oplus T^*M,[\cdot,\cdot]_H)$ be a Dirac structure
and let $\iota:S\hookrightarrow M$ be a maximal integral submanifold
for the (generalized) distribution $\Delta=\pi(L)\subset TM$.  Then
along $S$, we have $L=L(\Delta,\veps)$ for a unique
$\veps\in\Omega^2(S)$, and by the same argument as in
Proposition~\ref{regular}, we obtain
\[
\iota^*H=d\veps.
\]
Therefore we see that a Dirac structure induces a (generalized)
foliation of the manifold by trivializations $\LL=(\iota,\veps)$.

Note that in this example, $\tau_\LL = L|_S$ inherits a Lie
algebroid structure over $S$, since any sections $u,v\in
C^\infty(S,\tau_\LL)$ may be extended to $\tilde u,\tilde v\in
C^\infty(M,L)$ and then the expression
\[
[u,v]:=[\tilde u,\tilde v]|_S
\]
is independent of extension and defines a Lie bracket.  For general
trivializations, however, the ambient Dirac structure $L$ is
unavailable, and the argument fails.
\end{example}

A complex submanifold $S\subset M$ of a complex manifold is defined
by the property that $J(TS)=TS$.  Similarly, we define a
compatibility condition between a Courant trivialization and a
generalized complex structure.
\begin{defn}
A Courant trivialization  $\LL=(\iota,s)$ is said to be
\emph{compatible} with the generalized complex structure $\JJ$ if
and only if
\[
\JJ(\tau_\LL)=\tau_\LL,
\]
i.e. its generalized tangent bundle is a complex subbundle of $E$.
\end{defn}
An immediate consequence of the definition is that
$\pi(\JJ(N^*S))\subset TS$, which by Proposition~\ref{projpois} is
the statement that $P(N^*S)\subset TS$, i.e. $S$ is a coisotropic
submanifold for the Poisson structure $P$. Since $P$ is Poisson,
$\Delta=P(N^*S)$ integrates to a singular folation called the
\emph{characteristic foliation} of $S$.

Decomposing $\tau_\LL\otimes\CC$ into $\pm i$-eigenspaces for $\JJ$,
we obtain
\[
\tau_\LL\otimes\CC = \Ll\oplus\overline{\Ll}
\]
Note that the isotropic subbundle $\Ll\subset (E\otimes\CC)|_S$ is
contained in the ambient $+i$-eigenbundle $L$ of $\JJ$, i.e.
\[
\Ll\subset L|_S.
\]
Therefore, the argument of Example~\ref{extbra} concerning
restriction of Courant brackets applies and we obtain the following
result\footnote{This Lie algebroid was obtained independently by
Kapustin and Li~\cite{KapLi}, as defining the BRST complex
describing open strings with both ends on $\LL$.}.
\begin{prop}
Let $\JJ$ be a generalized complex structure and let $\LL$ be a
compatible Courant trivialization.  Define
$\Ll=\ker(\JJ-i)\cap(\tau_\LL\otimes\CC)$.  Then the Courant bracket
induces a Lie bracket on $C^\infty(S,\Ll)$, making
$(\Ll,[\cdot,\cdot],\pi)$ into a complex Lie algebroid over $S$.
\end{prop}
The associated Lie algebroid complex
$(C^\infty(S,\wedge^\bullet\Ll^*), d_\Ll)$  is actually
\emph{elliptic}, by the same reasoning as in
Proposition~\ref{ellipcx}, and may be used to study the deformation
theory of $\LL$, which we leave for a future work.

The Lie algebroid $\Ll$ projects to a generalized distribution
$A=\pi(\Ll)\subset TS\otimes\CC$, which is integrable and satisfies
$A+\overline A = TS\otimes\CC$.  The intersection $A\cap \overline A
= \Delta\otimes\CC$ coincides with the characteristic distribution
of the coisotropic submanifold $S$. Therefore, by the reasoning in
Proposition~\ref{transcx}, wherever $\Delta$ has constant rank, $A$
defines an invariant integrable holomorphic structure transverse to
the characteristic foliation.
\begin{corollary}\label{littlel}
Let $\LL$ be a compatible trivialization and let $\Ll$ be the
complex Lie algebroid defined above.  In a neighbourhood where the
characteristic distribution is of constant rank, $A=\pi(\Ll)\subset
TS\otimes\CC$ defines an integrable holomorphic structure transverse
to the characteristic foliation, which descends to the leaf space.
\end{corollary}
Since $\Ll$ is a Lie algebroid, we may associate to any compatible
trivialization $\LL$ the category of $\Ll$-modules, i.e. complex
vector bundles $V$ over $S$, equipped with flat Lie algebroid
connections with respect to $\Ll$.  We call these \emph{generalized
complex branes} on $\LL$.
\begin{defn}[Generalized complex brane] Let $\JJ$ be a generalized
complex structure and $\LL$ a compatible trivialization.  A
\emph{generalized complex brane} supported on $\LL$ is a module over
the Lie algebroid $\Ll$.
\end{defn}

We now indicate that there is a natural pullback map taking
generalized holomorphic bundles ($L$-modules) over $M$ to
$\Ll$-modules over $S$.  As a result, any compatible Courant
trivialization immediately supports not only the trivial brane
$V=\CC\times S$ but also the canonical brane $\iota^*K$, where $K$
is the canonical line bundle.
\begin{prop}
Let $\LL=(\iota,s)$ be a compatible trivialization and $V$ a
generalized holomorphic bundle over $M$.  Then $\iota^*V$ is
naturally a $\Ll$-module, and hence a generalized complex brane.
\end{prop}
\begin{proof}
Let $j:\Ll\hookrightarrow L|_S$ denote the inclusion, and let
$\delbar:C^\infty(V)\lra C^\infty(L^*\otimes V)$ be the flat Lie
algebroid connection defining the $L$-module structure. Then for
$v\in C^\infty(S,\iota^*V)$, choose an extension $\tilde v\in
C^\infty(M,V)$ and define $D:C^\infty(S,\iota^*V)\lra
C^\infty(S,\Ll^*\otimes\iota^*V)$ by
\[
Dv := j^*(\delbar\tilde v)|_S.
\]
This is independent of extension since $\pi(\Ll)\subset
TS\otimes\CC$, and is easily seen to be a flat $\Ll$-connection.
\end{proof}

We now describe the detailed structure of generalized complex branes
in the extremal cases of complex and symplectic geometry.
\begin{example}[Complex branes]
Let $\LL=(\iota,F)$, for $\iota:S\hookrightarrow M$ and $F\in
\Omega^2(S)$, be a generalized complex trivialization in a complex
manifold, so that $\iota^*H=dF$ and $\tau_\LL$, given
by~\eqref{exptau}, is a complex subbundle for
\begin{equation*}
\JJ_J=\left(\begin{matrix}-J&\\&J^*\end{matrix}\right).
\end{equation*}
This happenes if and only if
\begin{itemize}
\item $TS\subset TM$ is a complex subbundle for $J$, i.e. $S$ is a complex
submanifold, and
\item $J^*i_XF+i_{JX}F\in N^*S$ for all $X\in TS$, i.e. $F$ is
of type $(1,1)$.
\end{itemize}
In this case, the Lie algebroid $\Ll$ is given by
\[
\Ll = \{X+\xi\in T_{0,1}S \oplus T^*_{1,0}M\ :\ \iota^*\xi = i_X
F\},
\]
and is therefore isomorphic to $T_{0,1}S\oplus N^*_{1,0}S$, where
$N^*_{1,0}S$ denotes the holomorphic conormal bundle of $S$.  As a
result, a generalized complex brane supported on $\LL$ consists of a
holomorphic vector bundle $V$ over $S$, together with a holomorphic
section $\phi:V\lra N_{1,0}S\otimes V$ satisfying
\[
\phi\wedge\phi = 0  \in \wedge^2 N_{1,0}S\otimes \End(V).
\]
One sees directly from the description in Example~\ref{gcholb} of
generalized holomorphic bundles that $\Ll$-modules may be obtained
by pullback of $L$-modules.
\end{example}

\begin{example}[Symplectic branes]\label{sympbran}
As in the previous example, let $\LL=(\iota,F)$ be a compatible
trivialization, but for a symplectic structure
\begin{equation*}
\JJ_\omega=\left(\begin{matrix}&-\omega^{-1}\\\omega&\end{matrix}\right).
\end{equation*}
If $F=0$, then $\tau_\LL = TS\oplus N^*S$, and
$\JJ(\tau_\LL)=\tau_\LL$ is simply the requirement that
$\omega^{-1}(N^*S)\subset TS$ and $\omega(TS)\subset N^*S$, i.e. $S$
is a Lagrangian submanifold.  In this particular case,
$\Ll=\{X-i\omega(X)\ :\ X\in TS\}$, so that $\Ll$ is isomorphic as a
Lie algebroid to $TS$ itself;  hence $\Ll$-modules are simply flat
vector bundles supported on $S$.

However, there are symplectic branes beyond the flat bundles over
Lagrangians if we allow $F\neq 0$; in general, as we saw in
Corollary~\ref{littlel}, $S$ must be coisotropic and the Lie
algebroid $\Ll$ determines a complex distribution $A=\pi(\Ll)$
defining an invariant holomorphic structure transverse to the
characteristic foliation of $S$.  However for a symplectic
trivialization we have explicitly $\Ll = (\tau_\LL\otimes\CC)\cap
\Gamma_{-i\omega}$, and hence
\begin{equation}\label{ellexpress}
\Ll = \{X - i\omega(X)\in (TS\oplus T^*M)\otimes\CC\ :\
i_X(F+i\iota^*\omega)=0\}.
\end{equation}
Since $A=\pi(\Ll)$ and $\Delta\otimes\CC=A\cap\overline A$ defines
the characteristic foliation,~\eqref{ellexpress} implies that
$F+i\iota^*\omega$ is basic with respect to the foliation and
defines a closed, nondegenerate $(0,2)$-form on the leaf space.
Hence the leaf space inherits a natural holomorphic symplectic
structure.  In this way we obtain precisely the structure of
\emph{coisotropic A-brane}, discovered by Kapustin and
Orlov~\cite{Kapustin} in their search for geometric objects of the
Fukaya category beyond the well-known Lagrangian ones.

For such coisotropic trivializations, $\Ll$ is isomorphic as a Lie
algebroid to the distribution $A=\pi(\Ll)$, and so branes are vector
bundles equipped with flat partial $A$-connections. This implies
that they are flat along the characteristic distribution,
transversally holomorphic and invariant along the distribution.
Holomorphic bundles pulled back from the leaf space would provide
examples.
\end{example}

\begin{example}[Space-filling symplectic brane] A special case of
the preceding example is when the submanifold $S$ coincides with $M$
itself; then any brane over $\LL=(\text{id},s)$ is said to be
\emph{space-filling}. By the preceding argument, $\LL$ may be
described by a closed 2-form $F\in\Omega^2(M)$ such that
\[
\sigma = F+ i\omega
\]
defines a holomorphic symplectic structure on $M$, with complex
structure given by $J=-\omega^{-1}F$.

If $M$ supports such a space-filling brane, then any complex
submanifold $\iota:S\hookrightarrow M$ which is also coisotropic
with respect to $\sigma$ (for example, a complex hypersurface)
defines a compatible trivialization $\LL'=(\iota,\iota^*F)$ in $M$,
and we may produce examples of branes on $\LL'$ by pullback. The
holomorphic symplectic structure on its leaf space is also known as
the holomorphic symplectic reduction of $\LL'$.
\end{example}

\begin{example}[General space-filling branes]
The existence of a space-filling generalized complex brane places a
strong constraint on the generalized complex structure. Indeed, the
generalized tangent bundle $\tau_\LL$ determines an integrable
isotropic splitting of the Courant algebroid
\[
E = T^*\oplus\tau_\LL,
\]
so that the curvature $H$ vanishes.  If $\JJ$ is the generalized
complex structure, the constraint $\JJ(\tau_\LL)=\tau_\LL$ implies
that $\JJ$ must have upper triangular form in this splitting:
\[
\JJ = \begin{pmatrix}-J&P\\&J^*\end{pmatrix}.
\]
Here we use the canonical identification $\tau_\LL=TM$. Since $\JJ$
is upper triangular, $J$ is an integrable complex structure, for
which $T_{0,1} = \Ll$.  The real Poisson structure $P$ is of type
$(2,0)+(0,2)$ as can be seen from the fact $JP=PJ^*$, and the
complex bivector $\beta  = \tfrac{-1}{4}(Q+iP)$, for $Q=PJ^*$, is
such that the $+i$-eigenbundle of $\JJ$ can be written as
\[
L = T_{0,1}\oplus\Gamma_\beta,
\]
where $\beta$ is viewed as a map $T^*_{1,0}\lra T_{1,0}$. Courant
integrability then requires that $\beta$ be a holomorphic Poisson
structure.  Therefore we see that space-filling branes only exist
when $\JJ$ is a holomorphic Poisson deformation of a complex
manifold.
\end{example}

In fact, one can show using a combination of the above arguments, or
as is done in~\cite{zambon} by developing a theory of brane
reduction, that for an arbitrary generalized complex brane, the
Poisson and holomorphic structures transverse to the characteristic
foliation (when it is regular) are compatible, defining an invariant
transverse holomorphic Poisson structure.

Interesting relations between the coisotropic branes discussed in
this section and noncommutative geometry have appeared
in~\cite{KapA},\cite{Aldi} in particular;  for more on this
connection as well as the relation between coisotropic branes and
generalized K\"ahler geometry, see~\cite{GMN}.

\section{Appendix}\label{appendix}

\begin{prop}\label{restricte}
Let $E$ be an exact Courant algebroid over $M$ with \v{S}evera class
$[H]$, and suppose $\iota:S\hookrightarrow M$ is a submanifold. Then
\[
\iota^*E:= K^\bot/K,
\]
for $K=\Ann(TS)\subset E|_S$, inherits the structure of an exact
Courant algebroid over $S$ with \v{S}evera class $\iota^*[H]$.
\end{prop}
\begin{proof} We first show that $\iota^*E$ inherits a bracket.
Let $u,v\in C^\infty(S, K^\perp/K)$, and choose representatives
$u',v'\in C^\infty(S, K^\perp)$.  Extend these over $M$ as sections
$\tilde u, \tilde v\in C^\infty(M, E)$.  We claim that $[\tilde
u,\tilde v]|_S$ defines a section of $\iota^*E$ which is independent
of the choices made.

Firstly we observe that $[\tilde u,\tilde v]|_S\in
C^\infty(S,K^\perp)$, since $\pi[\tilde u,\tilde v] = [\pi\tilde
u,\pi \tilde v]$ and if $X,Y$ are vector fields tangent to $S$ then
$[X,Y]$ is also tangent to $S$.

Secondly we claim that $[\tilde u,\tilde v] + K$ is independent of
the choices made: for $p,q\in C^\infty(E)$ with $p|_S, q|_S \in
C^\infty(S,K)$, we have
\[
[\tilde u + p, \tilde v + q] - [\tilde u, \tilde v] = [\tilde u, q]
+ [p,\tilde v] + [p,q].
\]
Given any $x\in C^\infty(E)$ with $\pi(x)|_S\in C^\infty(S,TS)$, we
verify that $\IP{x,[\tilde u, q]} = \pi(\tilde u)\IP{x,q} -
\IP{[\tilde u,x],q}$ vanishes upon restriction to $S$, since
$\IP{x,q}$ vanishes along $S$ and $\pi([\tilde u,x])$ is tangent to
$S$. Similarly for the other two terms.  This shows that $[\tilde
u,\tilde v] + K$ is independent of the choices made.  The remainder
of the Courant algebroid properties are easily verified.
\end{proof}

\begin{prop}\label{pulldirac}
Let $L\subset E$ be a Dirac structure and assume that
\[
\iota^*L:= \frac{L\cap K^\perp + K}{K}
\]
is a smooth subbundle of $\iota^*E$.  Then it is a Dirac structure.
\end{prop}
\begin{proof}
We need only verify that the maximal isotropic subbundle $\iota^*L$
is involutive.  Let $u,v\in C^\infty(S, (L\cap K^\perp + K)/K)$. By
the definition of the Courant bracket in
Proposition~\ref{restricte}, we may choose representatives $u',v'\in
C^\infty(S,L\cap K^\perp + K)$ for $u,v$ and extend these as
sections of $E$ over $M$ in any way. In a neighbourhood $U\subset S$
where $L\cap K^\perp$ has constant rank, write $u' = x'+p'$, $v' =
y' + q'$, where $x',y'\in C^\infty(U, L\cap K^\perp)$ and $p',q'\in
C^\infty(U,K)$.  Then choose extensions $x,y\in C^\infty(V,L)$ for
$x',y'$ and $p,q\in C^\infty(V,E)$ for $p',q'$, where $V$ is an open
set in $M$ containing $U$.  Then
\[
[x+p,y+q] = [x,y] + [x,q] + [p,y] + [p,q].
\]
Since $x,y\in C^\infty(V,L)$ and $\pi(x),\pi(y)$ are tangent to $S$,
$[x,y]|_S\in C^\infty(U,L\cap K^\perp)$.  Also, $[x,q]|_S\in
C^\infty(U,K)$ since, for $z\in C^\infty(V,E)$ with $\pi(z)$ tangent
to $S$,
\[
\IP{z,[x,q]} = \pi(x)\IP{z,q} - \IP{[x,z],q},
\]
which vanishes along $S$ since $z$ and $[x,z]$ are both tangent to
$S$.  The same argument applies to show $[p,y]|_S, [p,q]|_S\in
C^\infty(S,K)$.  This proves that $[x+p,y+q]|_S\in C^\infty(S,L\cap
K^\perp + K)$, and hence that $\iota^*L$ is involutive in $U$. Since
$L\cap K^\perp$ has locally constant rank on an open dense set in
$S$, this argument shows that $\iota^*L$ is involutive, as required.
\end{proof}

\bibliographystyle{hplain}
\bibliography{gcg}

\end{document}